\documentclass[a4paper,12pt]{article}

\usepackage[font=small,labelfont=bf,width=0.8\textwidth,skip=16pt]{caption}
\usepackage{graphicx}
\usepackage{amsmath,amssymb}
\usepackage{here} 
\usepackage{bbm} 

\numberwithin{figure}{section}
\numberwithin{table}{section}
\numberwithin{equation}{section}

\renewcommand{\P}{\mathbb P}
\newcommand{\E}{\mathbb E}
\newcommand{\cF}{\mathcal F}
\newcommand{\cP}{\mathcal P}
\newcommand{\cL}{\mathcal L}
\newcommand{\cM}{\mathcal M}
\newcommand{\cC}{\mathcal C}

\newcommand{\I}{\mathbbm 1}
\newcommand{\N}{\mathbb N}
\newcommand{\Z}{\mathbb Z}
\newcommand{\R}{\mathbb R}

\newcommand{\bb}{\big|\big|}

\renewcommand\mid{\,\big|\,}

\newcommand{\bigtwomatrix}[4]{\begin{bmatrix}#1&#2\\#3&#4\end{bmatrix}}

\usepackage{amsthm}
\newtheoremstyle{mythm}{18pt}{0pt}{\itshape}{}{\bfseries}{.}{12pt}{}
\newtheoremstyle{mydefn}{18pt}{0pt}{}{}{\bfseries}{.}{12pt}{}
\theoremstyle{mythm}
\newtheorem{theorem}{Theorem}[section]
\newtheorem{lemma}[theorem]{Lemma}
\newtheorem{proposition}[theorem]{Proposition}
\newtheorem{conjecture}[theorem]{Conjecture}

\theoremstyle{mydefn}

\newtheorem{definition}[theorem]{Definition}
\newtheorem{example}[theorem]{Example}
\newtheorem{algorithm}[theorem]{Algorithm}

\begin{document}
\parindent 0cm
\parskip .5cm

{\centering
{\Huge Hidden Markov Models with\\[.5cm]
Multiple Observation Processes}\\[3cm]

{\large James Yuanjie Zhao}\\[4cm]

Submitted in total fulfilment of the requirements\\[-.2cm]
of the degree of Master of Philosophy\\[4cm]

Submitted 18 August 2010\\[-.2cm]
Revised 18 January 2011\\[.3cm]

Department of Mathematics and Statistics\\[-.2cm]
Department of Electrical and Electronic Engineering\\[-.2cm]
University of Melbourne\\[.3cm]

\vfill

}
\thispagestyle{empty}

\newpage
\section*{Abstract}
We consider a hidden Markov model with multiple observation processes, one of which is chosen at each point in time by a policy---a deterministic function of the information state---and attempt to determine which policy minimises the limiting expected entropy of the information state. Focusing on a special case, we prove analytically that the information state always converges in distribution, and derive a formula for the limiting entropy which can be used for calculations with high precision. Using this fomula, we find computationally that the optimal policy is always a threshold policy, allowing it to be easily found. We also find that the greedy policy is almost optimal.
\thispagestyle{empty}

\newpage
\section*{Declaration}
This is to certify that:
\vspace{-\parskip}
\begin{enumerate}
\item The thesis comprises only my original work towards the MPhil;
\item Due acknowledgement has been made to all other material used; and
\item The thesis is less than 50,000 words in length.
\end{enumerate}
\vspace{1cm}
\hspace{1.8cm}\includegraphics[width=.7\textwidth]{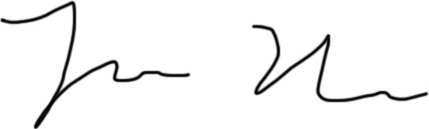}
\vspace{-1.8cm}

\begin{tabular}{c}
~\hspace{12.7cm}~\\
\hline
\end{tabular}
\vspace{1cm}
\thispagestyle{empty}

\newpage
\section*{Acknowledgements}
My deepest and sincerest appreciation goes to my supervisors, Bill Moran and Peter Taylor, for their countless hours of guidance, both in relation to this thesis and in more general matters.
\thispagestyle{empty}

\newpage
\tableofcontents
\thispagestyle{empty}

\newpage
\listoffigures
\thispagestyle{empty}

\newpage
\section{Introduction}
A hidden Markov model is an underlying Markov chain together with an imperfect observation on this chain. In the case of multiple observations, the classical model assumes that they can be observed simultaneously, and considers them as a single vector of observations. However, the case where not all the observations can be used at each point in time often arises in practical problems, and in this situation, one is faced with the challenge of choosing which observation to use.

We consider the case where the choice is made as a deterministic function of the previous information state, which is a sufficient statistic for the sequence of past observations. This function is called the \textit{policy}, which we rank according to the information entropy of the information state that arises due to that policy.

Our main results are:
\vspace{-\parskip}
\begin{itemize}
\item The information state converges in distribution for almost every underlying Markov chain, as long as each observation process gives a perfect information observation with positive probability;
\item In a special case (see Section \ref{sec:specialcase} for a precise definition), we can write down the limiting entropy explicitly as a rational function of subgeometric infinite series, which allows the calculation of limiting entropy to very good precision;
\item Computational results suggest that the optimal policy is a threshold policy, hence finding the optimal threshold policy is sufficient for finding the optimal policy in general;
\item Finding a locally optimal threshold policy is also sufficient, while finding a locally optimal general policy is sufficient with average probability 0.98; and
\item The greedy policy is optimal 96\% of the time, and close to optimal the remaining times, giving a very simple yet reasonably effective suboptimal alternative.
\end{itemize}

\newpage
\subsection{Motivation}
The theory of hidden Markov models was first introduced in a series of papers from 1966 by Leonard Baum and others under the more descriptive name of \textit{Probabilistic Functions of Markov Chains} \cite{baum}. An application of this theory was soon found in speech recognition, spurring development, and the three main problems---probability calculation, state estimation and parameter estimation---had essentially been solved by the time of Lawrence Rabiner's influential 1989 tutorial paper \cite{rabiner}.

The standard hidden Markov model consists of an underlying state which is described by a Markov chain, and an imperfect observation process which is a probabilistic function of this underlying state. In most practical examples, this single observation is equivalent to having multiple observations, since we can simply consider them as a single vector of simultaneous observations. However, this requires that these multiple observation can be made and processed simultaneously, which is often not the case.

Sometimes, physical constraints may prevent the simultaneous use of all of the available observations. This is most evident with a sensor which can operate in multiple modes. For example, a radar antenna must choose a waveform to transmit; each possible waveform results in a different distribution of observations, and only one waveform can be chosen for each pulse. Another example might be in studying animal populations, where a researcher must select locations for a limited pool of detection devices such as traps and cameras.

Even when simultaneous observations are physically possible, other constraints may restrict their availability. For example, in an application where processors are much more expensive than sensors, a sensor network might reasonably consist of a large number of sensors and insufficient processing power to analyse the data from every sensor, in which case the processor must choose a subset of sensors from which to receive data. Similarly, a system where multiple sensors share a limited communication channel must decide how to allocate bandwidth, in a situation where each bit of bandwidth can be considered a virtual sensor, not all of which can be simultaneously used.

Another example is the problem of searching for a target which moves according to a Markov chain, where observation processes represent possible sites to be searched. Indeed, MacPhee and Jordan's \cite{macphee} special case of this problem exactly corresponds to the special case we consider in Section \ref{sec:specialcase}, although with a very different cost function. Johnston and Krishnamurthy \cite{johnston} show that this search problem can be used to model file transfer over a fading channel, giving yet another application for an extended hidden Markov model with multiple observation processes.

Note that in the problem of choosing from multiple observation processes, it suffices to consider the case where only one observation is chosen, by considering an observation to be an allowable subset of sensors. The three main hidden Markov model problems of probability calculation, state estimation and parameter estimation remain essentially the same, as the standard algorithms can easily be adapted by replacing the parameters of the single observation process by those of whichever observation process is chosen at each point in time.

Thus, the main interesting problem in the hidden Markov model with multiple observation processes is that of determining the optimal choice of observation process, which cannot be adapted from the standard theory of hidden Markov models since it is a problem that does not exist in that framework. It is this problem which will be the focus of our work.

We will use information entropy of the information state as our measure of optimality. While Evans and Krishnamurthy \cite{evans} use a distance between the information state and the underlying state, it is not necessary to consider this underlying state explicitly, since the information state is by definition an unbiased estimator of the distribution of the underlying state. We choose entropy over other measures such as variance since it is a measure of uncertainty which requires no additional structure on the underlying set.

The choice of an infinite time horizon is made it order to simplify the problem, as is our decision to neglect sensor usage costs. These variables can be considered in future work.

\newpage
\subsection{Past Work}

The theory of hidden Markov models is already well-developed \cite{rabiner}. On the other hand, very little research has been done into the extended model with multiple observation processes. The mainly algorithmic solutions in the theory of hidden Markov models with a single observation process cannot be extended to our problem, since the choice of observation process does not exist in the unextended model.

Similarly, there is a significant amount of work on the sensor scheduling literature, but mostly considering autoregressive Gaussian processes such as in \cite{cochran}. The case of hidden Markov sensors was considered by Jamie Evans and Vikram Krishnamurthy in 2001 \cite{evans}, using policies where an observation process is picked as a deterministic function of the previous observation, and with a finite time horizon. They transformed the problem of choosing an observation into a control problem in terms of the information state, thereby entering the framework of stochastic control. They were able to write down the optimal policy as an intractible dynamic programming problem, and suggested the use of approximations to find the solution.

Krishnamurthy \cite{krishnamurthy} followed up this work by showing that this dynamic programming problem could be solved using the theory of Partially Observed Markov Decision Processes when the cost function is of the form
\[C(z)=\sum_iz(i)\big|\big|\delta(i)-z\big|\big|,\]
where $z$ is the information state, $\delta(i)\in\cP(S)$ is the Dirac measure and $||\cdot||$ is a piecewise constant norm. It was then shown that such piecewise linear cost functions could be used to approximate quadratic cost functions, in the sense that a sufficiently fine piecewise linear approximation must have the same optimal policy. In particular, this includes the Euclidean norm on the information state space, which corresponds to the expected mean-square distance between the information state and the distribution of the underlying chain. However, no bounds were found on how fine an approximation is needed.

The problem solved by Evans and Krishnamurthy is a similar but different problem to ours. We consider policies based on the information state, which we expect to perform better than policies based on only the previous observation, as the information state is a sufficient statistic for the sample path of observations (see Proposition \ref{prop:sufficientstatistic}, also \cite{witsenhausen}). We also consider and infinite time horizon, and specify information entropy of the information state as our cost function. Furthermore, while Evans and Krishnamurthy consider the primary tradeoff as that between the precision of the sensors and the cost of using them, we do not consider usage costs and only aim to minimise the uncertainty associated with the measurements.

Further work by Krishnamurthy and Djonin \cite{djonin} extended the set of allowable cost functions to a Lipschitz approximation to the entropy function, and proved that threshold policies are optimal under certain very restrictive assumptions. Their breakthrough uses lattice theory methods \cite{topkis} to show that the cost function must be monotonic in a certain way with respect to the information state, and thus the optimal choice of observation process must be characterised by a threshold. However, this work still does not solve our problem, as their cost function, a time-discounted infinite sum of expected costs, differs significantly from our limiting expected entropy, and furthermore their assumptions are difficult to verify in practice.

Another similar problem was also considered by Mohammad Rezaeian \cite{rezaeian}, who redefined the information state as the posterior distribution of the underlying chain given the sample path of observations up to the previous, as opposed to current, time instant, which allowed for a simplification in the recursive formula for the information state. Rezaeian also transformed the problem into a Markov Decision Process, but did not proceed further in his description.

The model for the special case we consider in Section \ref{sec:specialcase} is an instance of the problem of searching for a moving target, which was partially solved by MacPhee and Jordan \cite{macphee} with a very different cost function -- the expected cumulative sum of prescribed costs until the first certain observation. They proved that threshold policies are optimal for certain regions of parameter space by analysing the associated fractional linear transformations. Unfortunately, similar approaches have proved fruitless for our problem due to the highly non-algebraic nature of the entropy function.

Our problem as it appears here was first studied in unpublished work by Bill Moran and Sofia Suvorova, who conjectured that the optimal policy is always a threshold policy. More extensive work was done in \cite{honoursthesis}, where it was shown that the information state converges in distribution in the same special case that we consider in Section \ref{sec:specialcase}. It was also conjectured that threshold policies are optimal in this special case, although the argument provided was difficult to work into a full proof. However, \cite{honoursthesis} contains a mistake in the recurrence formula for the information state distribution, a corrected version of which appears as Lemma \ref{lem:informationdistributionrecurrence}. The main ideas of the convergence proof still work, and are presented in corrected and improved form in Section \ref{sec:convergence}.

\newpage
\section{Analytic Results}
\subsection{Definitions\label{sec:definitions}}

We begin by precisely defining the model we will use. In particular, we will make all our definitions within this section, in order to expediate referencing. For the reader's convenience, Table~\ref{tbl:definitions} at the end of this section lists the symbols we will use for our model.

For a sequence $X_0,X_1,\ldots$ and any non-negative integer $t\in\Z^+$, we will use the notation $X_{(t)}$ to represent the vector $\big(X_0,\ldots,X_t\big)$.

\begin{definition}
\label{def:markovchain}
A \textbf{Markov Chain} \cite{meyn} is a stochastic process $\big(X_t\big){}_{t\in\Z^+}$, such that for all times $t\ge s$, all states $x$ and all measurable sets $A$,
\[\P\big(X_t\in A\big|X_s=x,\cF_s\big)=\P\big(X_t\in A\big|X_s=x\big),\]
where $\cF_t$ denotes the canonical filtration. We will consistently use the symbol $X_t$ to refer to an underlying Markov chain, and $\pi_t=\P(X_t)$ to denote its distribution.

In the case of a time-homogeneous, finite state and discrete time Markov chain, this simplifies to a sequence of random variables $\big(X_t\big){}_{t\in\Z^+}$ taking values in a common finite state space $S=\{1,\ldots,n\}$, such that for all times $t\in\Z^+$, $X_{t+1}$ is conditionally independent of $X_{(t-1)}$ given $X_t$, and the distributions of $X_{t+1}$ given $X_t$ does not depend on $t$.

In this case, there exists $n\times n$ matrix $T$, called the \textbf{Transition Matrix}, such that for all $i,j\in S$ and $t\in\Z^+$,
\begin{equation}
T_{ij}=\P\big(X_{t+1}=j\mid X_t=i\big)=\P\big(X_{t+1}=j\mid X_t=i,X_{(t-1)}\big).
\end{equation}
Since we mainly consider Markov chains which are time-homogeneous and finite state, we will henceforth refer to them as Markov chains without the additional qualifiers.
\end{definition}

\begin{definition}
\label{def:observationprocess}
An \textbf{Observation Process} on the Markov chain $\big(X_t\big){}_{t\in\Z^+}$ is a sequence of random variables $\big(Y_t\big){}_{t\in\Z^+}$ given by $Y_t=c(X_t,W_t)$, where $c$ is a deterministic function and $\big(W_t\big)_{t\in\Z^+}$ is a sequence of independent and identically distributed random variables which is also independent of the Markov chain $\big(X_t\big){}_{t\in\Z^+}$ \cite{elliott}.

As before, we will only consider observation processes which take values in a finite set $V=\{1,\ldots,m\}$. Similarly to before, there exists an $m\times n$ matrix $M$, which we call the \textbf{Observation Matrix}, such that for all $i,j\in S$, $k\in V$ and $t\in\Z^+$,
\begin{align}
\begin{split}
&M_{jk}=\P\big(Y_t=k\mid X_t=j\big)=\P\big(Y_t=k\mid X_t=j,X_{(t-1)},Y_{(t-1)}\big);\\
&T_{ij}=\P\big(X_{t+1}=j\mid X_t=i\big)=\P\big(X_{t+1}=j\mid X_t=i,X_{(t-1)},Y_{(t)}\big).
\end{split}
\end{align}
\end{definition}

\begin{figure}[H]
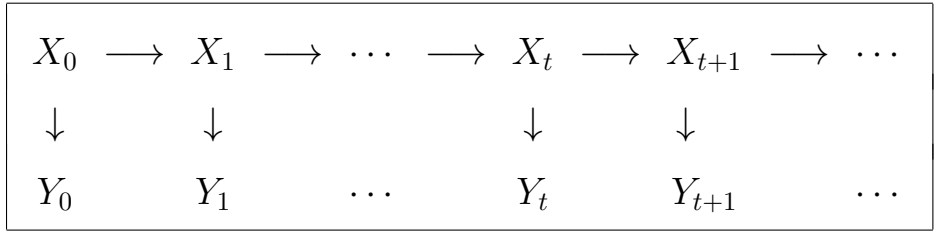

  \centering
  \large
  \begin{tabular}{|*{11}{@{~~}c}|}
    \hline
    $X_0$ & $\longrightarrow$ & $X_1$ & $\longrightarrow$ & $\cdots$ & $\longrightarrow$ & $X_t$ & $\longrightarrow$ & $X_{t+1}$ & $\longrightarrow$ & $\cdots$ \\[-6pt]
    $\downarrow$ && $\downarrow$ &&&& $\downarrow$ && $\downarrow$~~~ && \\[-6pt]
    $Y_0$ &                   & $Y_1$ &                   & $\cdots$ &                   & $Y_t$ &                   & $Y_{t+1}$ &                   & $\cdots$ \\
    \hline
  \end{tabular}
  \caption[A hidden Markov model]{An observation process $\big(Y_t\big)$ on a Markov chain $\big(X_t\big)$. At each node $X_t$, everything after $X_t$ is conditionally independent of everything before $X_t$, given $X_t$.}
  \label{fig:observationprocess}
\end{figure}

Heuristically, these two conditions can be seen as requiring that observations depend only on the current state, and do not affect future states. A diagrammatic interpretation is provided in Figure~\ref{fig:observationprocess}.

\vspace{3\parskip}
\vspace\parskip

Traditionally, a hidden Markov model is defined as the pair of a Markov chain and an observation process on that Markov chain. Since we will consider hidden Markov models with multiple observation processes, this definition does not suffice. We adjust it as follows.

\begin{definition}
\label{def:hiddenmarkovmodel}
A \textbf{Hidden Markov Model} is the triple of a Markov chain $\big(X_t\big){}_{t\in\Z^+}$, a finite collection of observation processes $\big\{\big(Y^{(i)}_t\big){}_{t\in\Z^+}\big\}_{i\in O}$ on $\big(X_t\big){}_{t\in\Z^+}$, and an additional sequence of random variables $\big(I_t\big){}_{t\in\Z^+}$, called the \textbf{Observation Index}, mapping into the index set $O$.

Note that this amends the standard definition of a hidden Markov model. For convenience, we will no longer explicitly specify our hidden Markov models to have multiple observation processes.
\end{definition}

It makes sense to think of $\big(X_t\big){}_{t\in\Z^+}$ as the state of a system under observation, $\big\{\big(Y^{(i)}_t\big){}_{t\in\Z^+}\big\}_{i\in O}$ as a collection of potential observations that can be made on this system, and $\big(I_t\big){}_{t\in\Z^+}$ as a choice of observation for each point in time.

Since our model permits only one observation to be made at each point in time, and we will wish to determine which one to use based on past observations, it makes sense to define $\big(I_t\big){}_{t\in\Z^+}$ as a sequence of random variables on the same probability space as the hidden Markov model.

We will discard the potential observations which are not used, leaving us with a single sequence of random variables representing the observations which are actually made.

\begin{definition}
\label{def:actualobservation}
The \textbf{Actual Observation} of a hidden Markov model is the sequence of random variables $\big(Y^{(I_t)}_t\big){}_{t\in\Z}$.

We will write $Y_t$ to mean $Y^{(I_t)}_t$, noting that this is consistent with our notation for a hidden Markov model with a single observation process $Y_t$. On the other hand, for a hidden Markov model with multiple observation processes, the actual observation $\big(Y_t\big){}_{t\in\Z^+}$ is not itself an observation process in general.
\end{definition}

Since our goal is to analyse a situation in which only one observation can be made at each point in time, we will consider our hidden Markov model as consisting only of the underlying state $\big(X_t\big){}_{t\in\Z^+}$ and the actual observation $\big(Y_t\big){}_{t\in\Z^+}$. Where convenient, we will use the abbreviated terms \textit{state} and \textit{observation} at time $t$ to mean $X_t$ and $Y_t$ respectively.

For any practical application of this model to a physical system, the underlying state cannot be determined, otherwise there would be no need to take non-deterministic observations. Therefore, we need a way of estimating the underlying state from the observations.

\begin{definition}
\label{def:informationstate}
The \textbf{Information State Realisation} of a hidden Markov model at time $t$ is the posterior distribution of $X_t$ given the actual observations and observation indices up to time $t$.

To make this definition more precise, we introduce some additional notation.

First, recall that $\big(X_t\big){}_{t\in\Z^+}$ has state space $S=\{1,\ldots,n\}$, and define the set of probability measures on $S$,
\begin{equation}
\cP(S)\cong\Big\{(p_1,\ldots,p_n)\in\R^n:p_i\ge0\;\forall\;i\in S,\;\textstyle\sum_{i\in S}p_i=1\Big\}.
\end{equation}
Second, for a random variable $X$ with state space $S$, and an event $E$, define the posterior distribution of $X$ given $E$,
\begin{equation}
\P\big(X\mid E\big)=\Big(\P\big(X=1\mid E\big)\;,\;\ldots\;,\;\P\big(X=n\mid E\big)\Big)\in\cP(S).
\end{equation}
Although we make this definition in general, we purposely choose the letters $X$ and $S$, coinciding with the letters used to represent the underlying Markov chain and the state space, as this is the context in which we will use this definition. Then, the information state realisation is a function
\begin{align}
\begin{split}
&\hspace*{1.33cm}z_t:V^{t+1}\times O^{t+1}\longrightarrow\cP(S),\\
&z_t\big(y_{(t)};i_{(t)}\big)=\P\big(X_t\mid Y_{(t)}=y_{(t)},I_{(t)}=i_{(t)}\big).
\end{split}
\end{align}
\end{definition}

This extends very naturally to a random variable.

\begin{definition}
\label{def:informationrv}
The \textbf{Information State Random Variable} is
\begin{equation*}
Z_t=z_t\big(Y_{(t)};I_{(t)}\big)=\P\big(X_t\mid Y_{(t)},I_{(t)}\big).
\end{equation*}
Its distribution is the \textbf{Information State Distribution} $\mu_t=\P(Z_t)$, taking values in $\cP(\cP(S))$, the space of Radon probability measures on $\cP(S)$, which is a subset of the real Banach space of signed Radon measures on $\cP(S)$.

Thus, the information state realisation is exactly a realisation of the information state random variable. It is useful because it represents the maximal information we can deduce about the underlying state from the observation index and the actual observation, as shown in Proposition \ref{prop:sufficientstatistic}.

For the purpose of succinctness, we will refer to any of $z_t$, $Z_t$ and $\mu_t$ as simply the \textbf{Information State} when the context is clear.
\end{definition}

\begin{definition}
\label{def:sufficientstatistic}
A random variable $Z$ is a \textbf{sufficient statistic} for a parameter $X$ given data $Y$ if for any values $y$ and $z$ of $Y$ and $Z$ respectively, the probability $\P\big(Y=y\mid Z=z,X=x\big)$ is independent of $x$ \cite{statsbook}. As before, we make the definition in general, but purposely choose the symbols $X$, $Y$ and $Z$ to coincide with symbols already defined.

In our case, $X$, which is a random variable, is used in the context of a parameter. Our problem takes place in a Bayesian framework, where the information state represents our belief about the underlying state, and is updated at each observation.
\end{definition}

\begin{proposition}
\label{prop:sufficientstatistic}
~The information random variable $Z_t$ is a sufficient statistic for the underlying state $X_t$, given the actual observations $Y_{(t)}$ and the observation indices $I_{(t)}$.
\end{proposition}
\begin{proof}
By Definition \ref{def:sufficientstatistic}, we need to prove that for all $y\in V^{t+1}$ and $i\in O^{t+1}$,
\begin{equation}
\P\big(Y,I\mid Z,x\big)=\P\big(Y_{(t)}=y,I_{(t)}=i\mid Z_t=z_t(y\,;i),X_t=x\big)
\end{equation}
is independent of $x$.

First, note that the event $\{Z_t=z_t(y\,;i)\}$ is the disjoint union of events $\{Y_{(t)}=y',I_{(t)}=i'\}$ over all $(y',i')\in V^{t+1}\times O^{t+1}$ such that $z_t(y'\,;i')=z_t(y\,;i)$.

Next, if $z_t(y\,; i)=z_t(y'\,; i')$, then by definition of $z_t$, for all $x\in S$,
\begin{equation}
\P\big(X_t=x\mid Y_{(t)}=y,I_{(t)}=i\big)=\P\big(X_t=x\mid Y_{(t)}=y',I_{(t)}=i'\big).
\end{equation}
Then, by definition of conditional probability,
\begin{equation}
\frac{\P\big(X_t=x,Y_{(t)}=y',I_{(t)}=i'\big)}{\P\big(X_t=x,Y_{(t)}=y,I_{(t)}=i\big)}=\frac{\P\big(Y_{(t)}=y',I_{(t)}=i'\big)}{\P\big(Y_{(t)}=y,I_{(t)}=i\big)}.
\end{equation}
Hence,
\begin{align}
\P\big(Y,I\mid Z,x\big)&=\P\big(Y_{(t)}=y,I_{(t)}=i\mid Z_t=z_t(y\,;i),X_t=x\big)\nonumber\\
&=\frac{\P\big(Y_{(t)}=y,I_{(t)}=i,Z_t=z_t(y\,;i),X_t=x\big)}{\P\big(Z_t=z_t\big(y\,;i\big),X_t=x\big)}\nonumber\\[4pt]
&=\frac{\P\big(Y_{(t)}=y,I_{(t)}=i,X_t=x\big)}{\sum_{y',i'}\P\big(Y_{(t)}=y',I_{(t)}=i',X_t=x\big)}\nonumber\\[4pt]
&=\Bigg(\sum_{y',i'}\frac{\P\big(Y_{(t)}=y',I_{(t)}=i',X_t=x\big)}{\P\big(Y_{(t)}=y,I_{(t)}=i,X_t=x\big)}\Bigg)^{-1}\nonumber\\[4pt]
&=\Bigg(\sum_{y',i'}\frac{\P\big(Y_{(t)}=y',I_{(t)}=i'\big)}{\P\big(Y_{(t)}=y,I_{(t)}=i\big)}\Bigg)^{-1}.
\end{align}
Each sum above is taken over all $(y',i')\in V^{t+1}\times O^{t+1}$ such that $z_t(y'\,;i')=z_t(y\,;i)$. This expression is clearly independent of $x$, which completes the proof that $Z_t$ is a sufficient statistic for $X_t$.
\end{proof}

Since the information state represents all information that can be deduced from the past, it makes sense to use it to determine which observation process to use in future.

\begin{definition}
\label{def:policy}
A \textbf{policy} on a hidden Markov model is a deterministic function $g:\cP(S)\rightarrow O$, such that for all $t\in\Z^+$, $I_{t+1}=g(Z_t)$. We will use the symbol $A_i=g^{-1}\{i\}$ to denote the preimage of the observation method $i$ under the policy, that is, the subset of $\cP(\cP(S))$ on which observation method $i$ is prescribed by the policy. We will always consider the policy $g$ as fixed.
\end{definition}

Since $Z_t$ is a function of $Y_{(t)}$ and $I_{(t)}$, this means that $I_{(t+1)}$ is a function of $Y_{(t)}$ and $I_{(t)}$. Then by induction, we see that $I_{(t+1)}$ is a function of $Y_{(t)}$ and $I_0$. Therefore, if we prescribe some fixed $I_0$, then $I_t$ is a function of $Y_{(t)}$.

For fixed $I_0$, we can write
\begin{equation}
\label{eqn:noIdependence}
Z_t=\P\big(X_t\mid I_{(t)},Y_{(t)}\big)=\P\big(X_t\mid Y_{(t)}\big).
\end{equation}
Hence, the information random variable is a deterministic function of only $Y_{(t)}$. In particuar, the information state $z_t$ can be written with only one argument, that is, $Z_t=z_t(Y_{(t)})$.

Since our aim is to determine the underlying state with the least possible uncertainty, we need to introduce a quantifier of uncertainty. There are many possible choices, especially if the state space has additional structure. For example, variance would be a good candidate in an application where the state space embeds naturally into a real vector space.

However, in the general case, there is no particular reason to suppose our state space has any structure; our only assumption is that it is finite, in which case information entropy is the most sensible choice, being a natural, axiomatically-defined quantifier of uncertainty for a distribution on a countable set without any additional structure \cite{cover}.

\begin{definition}
The \textbf{Information Entropy} of a discrete probability measure $(p_1,\ldots,p_n)\in\cP(S)$ is given by
\[H\big((p_1,\ldots,p_n)\big)=-\sum_jp_j\log p_j.\]
We will use the natural logarithm, and define $0\log0=0$ in accordance with the fact that $p\log p\rightarrow0$ as $p\rightarrow0$.
\end{definition}

Since $Z_t$ takes values in $\cP(S)$, $H(Z_t)$ is well-defined, and by definition measures the uncertainty in $X_t$ given $Y_{(t)}$, and therefore by Proposition \ref{prop:sufficientstatistic}, measures the uncertainty in our best estimate of $X_t$. Thus, the problem of minimising uncertainty becomes quantified as one of minimising $H(Z_t)$.

We are particularly interested in the limiting behaviour, and thus, the main questions we will ask are:
\begin{itemize}
\item Under what conditions, and in particular what policies, does $H(Z_t)$ converge as $t\rightarrow\infty$?
\item Among the policies under which $H(Z_t)$ converges, which policy gives the minimal limiting value of $H(Z_t)$?
\item Are there interesting cases where $H(Z_t)$ does not converge, and if so, can we generalise the above results?
\end{itemize}

\newpage
\begin{table}[H]
  \centering
  \begin{tabular}{|c|l|l|}
    \hline
    \textbf{Symbol}   & \textbf{Value}    & \textbf{Meaning} \\
    \hline
    $\cP(S)$      & set        & probability measures on $S$\\
    $\cP(\cP(S))$ & set & probability measures on $\cP(S)$\\
    $A_i$         & set & region of observation process $i$\\
    $H$           & function                & information entropy\\
    $I_t$         & random variable         & observation index\\
    $O$           & finite set              & set of observation processes\\
    $S$           & finite set              & state space of Markov chain\\
    $V$           & finite set              & observation space\\
    $W_t$         & random variable         & observation randomness\\
    $X_t$         & random variable         & Markov chain\\
    $~~Y^{(i)}_t$ & random variable         & observation process\\
    $Y_t$         & random variable         & actual observation $Y^{(I_t)}_t$\\
    $Z_t$         & random variable         & information state random variable\\
    $g$           & function                & policy\\
    $m$           & integer                 & number of observation values\\
    $n$           & integer                 & number of states\\
    $t$           & integer                 & position in time\\
    $z_t$         & distribution            & information state realisation\\
    $\pi_t$       & distribution            & Markov chain distribution\\
    $\mu _t$      & distribution            & information state distribution\\
    \hline
  \end{tabular}
  \caption{List of symbols, ordered alphabetically.}
  \label{tbl:definitions}
\end{table}

\newpage
\subsection{Convergence}
\label{sec:convergence}

In this section, we will prove that under certain conditions, the information state converges in distribution. This fact is already known for classical hidden Markov models, and is quite robust: LeGland and Mevel \cite{legland} prove geometric ergodicity of the information state even when calculated from incorrectly specified parameters, while Capp\'e, Moulines and Ryd\'en \cite{cappe} prove Harris recurrence of the information state for certain uncountable state underlying chains. We will present a mostly elementary proof of convergence in the case of multiple observation processes.

To determine the limiting behaviour of the information state, we begin by finding an explicit form for its one-step time evolution.

\begin{definition}
\label{def:rfunction}
For each observation process $i$ and each observed state $y$, the \textbf{$r$-function} is the function $r_{i,y}:\cP(S)\rightarrow\cP(S)$ given by
\begin{equation*}
r_{i,y}(z)=\frac{\sum_{x,j\in S}M^{(i)}_{x,y}T_{j,x}z_j\delta(x)}{\sum_{x,j\in S}M^{(i)}_{x,y}T_{j,x}z_j},
\end{equation*}
where $\delta:S\rightarrow\cP(S)$ is the Dirac measure on $S$ and $z_j$ is the $j$th component of $z\in\cP(S)\subset\R^n$.
\end{definition}

\begin{lemma}
\label{lem:informationstaterecurrence}
In a hidden Markov model with multiple observation processes and a fixed policy $g$, the information state satisfies the recurrence relation
\begin{equation}
z_{t+1}(y_{(t+1)})=r_{g(z_t(y_{(t)})),y_{t+1}}(z_t(y_{(t)})).
\end{equation}
\end{lemma}
\begin{proof}
Let $i_{t+1}=g(z_t(y_{(t)}))$ and $k_t=\P\big(Y^{(i_{(t)})}_{(t)}=y_{(t)}\big)$. By the Markov property as in Definition \ref{def:observationprocess}, and the simplification (\ref{eqn:noIdependence}),
\begin{align}
z_{t+1}(y_{(t+1)})_x
&=\P\big(X_{t+1}=x\mid Y^{(i_{(t+1)})}_{(t+1)}=y_{(t+1)}\big)\nonumber\\
&=\frac1{\P\big(Y^{(i_{(t+1)})}_{(t+1)}=y_{(t+1)}\big)} \sum_j\P\big(X_{t+1}=x,Y^{(i_{(t+1)})}_{(t+1)}=y_{(t+1)},X_t=j\big)\nonumber\\
&=\frac1{k_{t+1}}\sum_j\P\big(Y^{(i_{t+1})}_{t+1}=y_{t+1}\mid X_{t+1}=x\big)\nonumber\\[-.55cm]
&\hspace{3cm}\times\P\big(X_{t+1}=x\mid X_t=j\big)\P\big(X_t=j,Y^{(i_{(t)})}_{(t)}=y_{(t)}\big)\nonumber\\
&=\frac1{k_{t+1}}\sum_jM^{(i_{t+1})}_{x,y_{t+1}}T_{j,x}k_tz_t(y_{(t)})_j\nonumber\\
&=\frac{k_t}{k_{t+1}}\sum_jM^{(i_{t+1})}_{x,y_{t+1}}T_{j,x}z_t(y_{(t)})_j\nonumber\\
&=\frac{\sum_jM^{(i_{t+1})}_{x,y_{t+1}}T_{j,x}z_t(y_{(t)})_j}{\sum_x\sum_jM^{(i_{t+1})}_{x,y_{t+1}}T_{j,x}z_t(y_{(t)})_j},
\end{align}
since $k_t/k_{t+1}$ does not depend on $x$ and $\sum_xz_{t+1}(y_{(t+1)})_x=1$.
\end{proof}

Note that for each information state $z$ and each observation process $i$, there are at most $m$ possible information states at the next step, which are given explicitly by $r_{i,y}(z)$ for each observation $y\in V$.

\begin{lemma}
\label{lem:informationdistributionrecurrence}
The information distribution satisfies the recurrence relation
\begin{equation*}
\mu_{t+1}=\sum_{i\in O}\sum_{y\in V}\int_{A_i}\big(z\cdot T\cdot M^{(i)}\big)\!_y\,\delta\big(r_{i,y}(z)\big)d\mu_t(z),
\end{equation*}
where the sum is taken over all observation processes $i$ and all observation states $y$, $\delta:\cP(S)\rightarrow\cP(\cP(S))$ is the Dirac measure on $\cP(S)$, and $\cdot$ is the matrix product considering $z\in\cP(S)\subset\R^n$ as a row vector.
\end{lemma}
\begin{proof}
Since $Z_t=\cP\big(X_t\big|Y_{(t)}\big)$ is a deterministic function of $Y_{(t)}$, given that $Y_{(t+1)}=y_{(t+1)}$,
\begin{equation}
Z_{t+1}=z_{t+1}(y_{(t+1)})=r_{g(z_t(y_{(t)})),y_{t+1}}(z_t(y_{(t)})).
\end{equation}
This depends only on $z_t(y_{(t)})$ and $y_{t+1}$, so given that $Z_t=z$ and $Y_{t+1}=y$,
\begin{equation}
Z_{t+1}=r_{g(z),y}(z).
\end{equation}
Integration over $(Z_t,Y_{t+1})\in\cP(S)\times V$ gives
\begin{align}
\mu_{t+1}
&=\int_{\cP(S)}\sum_{y\in V}\delta\big(r_{g(z),y}(z)\big)\P\big(Y_{t+1}=y\big|Z_t=z\big)d\mu_t(z)\nonumber\\
&=\sum_{i\in O}\sum_{y\in V}\int_{A_i}\delta\big(r_{g(z),y}(z)\big)\P\big(Y_{t+1}=y\big|Z_t=z\big)d\mu_t(z).\label{eqn:mutplusone}
\end{align}
By Definition \ref{def:informationstate}, $Z_t$ is the posterior distribution of $X_t$ given the observations up to time $t$, so $\P\big(X_t=x\big|Z_t=z\big)=z_x$, the $x$th coordinate of the vector $z\in\cP(S)\subset\R^n$. Since $Z_t$ is a function of $Y_{(t)}$, which is a function of $X_{(t)}$ and the observation randomness $W_{(t)}$, by the Markov property as in Definition \ref{def:observationprocess},
\begin{align}
\P\big(Y_{t+1}=y\big|Z_t=z\big)
&=\sum_{x\in S}\big(Y_{t+1}=y\big|X_t=x\big)\P\big(X_t=x\big|Z_t=z\big)\nonumber\\
&=\sum_{x\in S}\big(T\cdot M^{(i)}\big)\!{}_{x,y}z_x=\big(z\cdot T\cdot M^{(i)}\big)\!{}_y.\label{eqn:pygivenz}
\end{align}
Substituting (\ref{eqn:pygivenz}) into (\ref{eqn:mutplusone}) completes the proof.
\end{proof}

Note that Lemma \ref{lem:informationdistributionrecurrence} shows that the information distribution is given by a linear dynamical system on $\cP(\cP(S))$, and therefore the information state is a Markov chain with state space $\cP(S)$. We will use tools in Markov chain theory to analyse the convergence of the information state, for which it will be convenient to give a name to this recurrence.

\begin{definition}
\label{def:transitionfunction}
The \textbf{transition function} of the information distribution is the deterministic function $F:\cP(\cP(S))\rightarrow\cP(\cP(S))$ given by $F(\mu_t)=\mu_{t+1}$, extended linearly to all of $\cP(\cP(S))$ by the recurrence in Lemma \ref{lem:informationdistributionrecurrence}. The coefficients $\alpha_{i,y}(z)=\big(z\cdot T\cdot M^{(i)}\big)\!{}_y$ are called the \textbf{$\alpha$-functions}.
\end{definition}

We now give a criterion under which the information state is always positive recurrent.

\begin{definition}
A discrete state Markov Chain $X_t$ is called \textbf{Ergodic} if it is irreducible, aperiodic and positive recurrent. Such a chain has a unique invariant measure $\pi$, which is a limiting distribution in the sense that $X_t$ converges to $\pi$ in total variation norm \cite{meyn}.
\end{definition}

\begin{definition}
A discrete state Markov Chain $X_t$ is called \textbf{Positive} if every transition probability is strictly positive, that is, for all $i,j\in S$, $\P(X_{t+1}=i|X_t=j)>0$. This is a stronger condition than ergodicity.
\end{definition}

\begin{definition}
\label{def:anchored}
We shall call a hidden Markov model \textbf{Anchored} if the underlying Markov chain $X_t$ is ergodic, and for each observation process $i$, there is a state $x_i$ and an observation $y_i$ such that $M^{(i)}_{x_i,y_i}>0$ and $M^{(i)}_{x,y_i}=0$ for all $x\ne x_i$. The pair $(x_i,y_i)$ is called an \textbf{Anchor Pair}.
\end{definition}

Heuristically, the latter condition allows for perfect information $\delta(x_i)$ whenever the observation $y_i$ is made using observation process $i$. This anchors the information chain in the sense that this state can be reached with positive probability from any other state, thus resulting in a recurrent atom in the uncountable state chain $Z_t$. On the other hand, since each information state can make a transition to only finitely many other information states, starting the chain at $\delta(x_i)$ results in a discrete state Markov chain, for which it is much easier to prove positive recurrence.

\begin{lemma}
\label{lem:anchor}
In an anchored hidden Markov model, for any anchor pair $(x_i,y_i)$, $r_{i,y_i}(z)=\delta(x_i)$ for all $z\in\cP(S)$.
\end{lemma}
\begin{proof}
When $x\ne x_i$, $M^{(i)}_{x,y_i}=0$ by Definition \ref{def:anchored}, so every term in the numerator of Definition \ref{def:rfunction} is zero except the coefficient of $\delta(x_i)$. Since we know the coefficients have sum 1, it follows that $r_{i,y_i}=\delta(x_i)$.
\end{proof}

\begin{lemma}
\label{lem:alpha}
In a positive, anchored hidden Markov model, the $\alpha$-functions $\alpha_{i,y_i}$, for each $i\in O$, are uniformly bounded below by some $\epsilon>0$, that is, $\alpha_{i,y_i}(z)\ge\epsilon$ for all $i$ and $z$.
\end{lemma}
\begin{proof}
We can write $\alpha_{i,y_i}(z)=\sum_xz_xT_{x,x_i}M^{(i)}_{x_i,y_i}$ by Definitions \ref{def:transitionfunction} and \ref{def:anchored}, which is bounded below by $\min_xT_{x,x_i}M^{(i)}_{x_i,y_i}$ since $\sum_xz_x=1$. Since each $M^{(i)}_{x_i,y_i}>0$, if all the entries of $T$ are positive, then $\alpha_{i,y_i}(z)$ is bounded below uniformly in $z$ for fixed $i$, which then implies a uniform bound in $z$ and $i$ since there are only finitely many $i$.
\end{proof}

\begin{definition}
\label{def:orbit}
For each state $x\in S$, the \textbf{Orbit} $R_x$ of $\delta(x)\in\cP(S)$ under the $r$-functions is
\begin{align}
R_x=\big\{\delta(x)\big\}&\cup\big\{r_{i,y}(\delta(x)):\alpha_{i,y}(\delta(z))>0\big\}\nonumber\\
&\cup\big\{r_{i',y'}\circ r_{i,y}(\delta(x)):\alpha_{i',y'}(r_{i,y}(\delta(x)))\alpha_{i,y}(\delta(x))>0\big\}\nonumber\\
&\cup\;\cdots.\nonumber
\end{align}
By requiring the $\alpha$-functions to be positive, we exclude points in the orbit which are reached with zero probability. Let $R=\bigcup_xR_x$.
\end{definition}

\begin{proposition}
\label{prop:discreteness}
In a positive, anchored hidden Markov model, there exists a constant $0<\lambda<1$ such that for all measures $Z\in\cP(\cP(S))$, the mass of the measure $F^t(Z)$ outside $R$ is bounded by $\lambda^t$, that is, $F^t(Z)(R^c)\le\lambda^t$.
\end{proposition}
\begin{proof}
We can rewrite Definition \ref{def:transitionfunction} as
\begin{equation}
F(Z)=\int_{\cP(S)}Q(z)dZ(z),
\end{equation}
where
\begin{equation}
Q(z)=\sum_i\sum_y\mathbbm1_{A_i}(z)\alpha_{i,y}(z)\delta\big(r_{i,y}(z)\big).
\end{equation}
In this notation, the integral is the Lebesgue integral of the function $Q$ with respect to the measure $Z$. Since $Q$ takes values in the $\cP(\cP(S))$ and $Z$ is a probability, the integral also takes values in $\cM(\cP(S))$, thus $F$ maps the information state space $\cP(\cP(S))$ to itself.

Since $Q(z)$ is a measure supported on the set of points reachable from $z$ via an $r$-function, and $R$ is a union of orbits of $r$-functions and therefore closed under $r$-functions, it follows that all mass in $R$ is mapped back into $R$ under the evolution function, that is
\begin{equation}
\bigg(\int_RQ(z)dZ(z)\bigg)(R)=Z(R)=1-Z(R^c).
\end{equation}
On the other hand, by Lemma \ref{lem:anchor}, $r_{i,y_i}(z)=\delta(x_i)\in R$ for all $z$, hence
\begin{align}
\bigg(\int_{R^c}Q(z)dZ(z)\bigg)(R)
&\ge\inf_z\big(Q(z)\big)(R)\nonumber\\
&\ge\inf_z\Big(\I_{A_{g(z)}}(z)\alpha_{g(z),y_{g(z)}}(z)\delta(r_{g(z),y_{g(z)}}(z))\Big)(R)\nonumber\\
&=\inf_z\alpha_{g(z),y_{g(z)}}(z)\ge\inf_{i,z}\alpha_{i,y_i}(z).
\end{align}
Putting these together gives
\begin{align}
F(Z)(R)
&=\bigg(\int_RQ(z)dZ(z)+\int_{R^c}Q(z)dZ(z)\bigg)(R)\nonumber\\
&\ge1-\big(1-\inf_{i,z}\alpha_{i,y_i}(z)\big)Z(R^c).
\end{align}
Setting $\lambda=1-\inf_{i,z}\alpha_{i,y_i}(z)$ gives $F(Z)(R^c)\le\lambda Z(R^c)$, hence $F^t(Z)(R^c)\le\lambda^t$ by induction. By Lemma \ref{lem:alpha}, $\lambda<1$, while $\lambda>0$ since we can always choose a larger value.
\end{proof}

Up to this point, we have considered the evolution function as a deterministic function $F:\cP(\cP(S))\rightarrow\cP(\cP(S))$. However, we can also consider it as a probabilistic function $F:\cP(S)\rightarrow\cP(S)$. By Definition \ref{def:transitionfunction}, $F$ maps points in $R$ to $R$, hence the restriction $F|_R:R\rightarrow R$ gives a probabilistic function, and therefore a Markov chain, with countable state space $R$.

By Proposition \ref{prop:discreteness}, the limiting behaviour of the information chain takes place almost entirely in $R$ in some sense, so we would expect that convergence of the restricted information chain $F|_R$ is sufficient for convergence of the full information chain $F$. This is proved below.

\begin{proposition}
\label{prop:positiverecurrent}
In a positive, anchored hidden Markov model, under any policy, the chain $F|_R$ has at least one state of the form $\delta(x_i)$ which is positive recurrent, that is, whose expected return time is finite.
\end{proposition}
\begin{proof}
Construct a Markov chain $P$ on the set $\{A_i\}\cup\{R_x\}$, with transition probabilities $P(A_i,R_{x_i})>0$ for all $i$, $P(R_x,A_i)>0$ whenever $R_x\cap A_i$ is nonempty, and all other transition probabilities zero. We note that this is possible since we allow each state a positive probability transition to some other state.

Since $P$ is a finite state Markov chain, it must have a recurrent state. Each state $R_x$ can reach some state $A_i$, so some state $A_i$ is recurrent; call it $A_1$.

Consider a state $z_0=r^{k_2}(\delta(x_2))\in R_{x_2}\subseteq R$ of the chain $F$ which is reachable from $\delta(x_1)$, where $r^{k_2}$ is a composition of $k_2$ $r$-functions with corresponding $\alpha$-functions nonzero. Since the $A_i$ partition $\cP(S)$, one of them must contain $z_0$; call it $A_3$. We will assume $A_3\ne A_1$; the proof follows the same argument and is simpler in the case when $z_0\in A_1$.

By definition of the $\alpha$-functions,
\begin{equation}
\P\big(Z_1=\delta(x_3)\big|Z_0=z_0\big)=\alpha_{3,y_3}(z_0)>0.
\end{equation}
This means that $\delta(x_3)$ is reachable from $\delta(x_1)$ in the chain $F$, hence in the chain $P$, $A_3$ is reachable from $A_1$, and by recurrence of $A_1$, $A_1$ must also be reachable from $A_3$ via some sequence of positive probability transitions
\begin{equation}
\label{eqn:path}
A_3\rightarrow R_{x_3}\rightarrow A_4\rightarrow R_{x_4}\rightarrow\cdots\rightarrow A_1.
\end{equation}
By definition of $P(R_{x_3},A_4)>0$, $R_{x_3}\cap A_4$ is nonempty, and thus contains some point $r^{k_3}(\delta(x_3))$, where $r^{k_3}$ is a composition of $k_3$ $r$-functions with corresponding $\alpha$ nonzero. 

By Definition \ref{def:orbit}, each transition $r^\ell(\delta(x_3))$ to $r^{\ell+1}(\delta(x_3))$ in the information chain occurs with positive probability, so
\begin{equation}
\P\big(Z_{k_3+1}=r^{k_3}(\delta(x_3))\big|Z_1=\delta(x_3)\big)=\beta_3>0.
\end{equation}
Since $r^{k_3}(\delta(x_3))\in A_4$, by anchoredness and positivity,
\begin{equation}
\P\big(Z_{k_3+2}=\delta(x_4)\big|Z_{k_3+1}=r^{k_3}(\delta(x_3))\big)=\gamma_3>0.
\end{equation}
The Markov property then gives
\begin{equation}
\P\big(Z_{k_3+2}=\delta(x_4)\big|Z_0=z_0\big)=\alpha_{3,y_3}(z_0)\beta_3\gamma_3>0.
\end{equation}
Continuing via the sequence (\ref{eqn:path}), we obtain
\begin{equation}
\P\big(Z_{k_1+\cdots+k_j+j+1}=\delta(x_1)\big|Z_0=z_0\big)=\alpha_{3,y_3}(z_0)\beta_3\gamma_3\cdots\beta_j\gamma_j>0.
\end{equation}
Thus, for every state $z\in R$ reachable from $\delta(x_1)$, we have found constants $s_z\in\N$ and $c_z>0$ such that
\[\P\big(Z_{s_z}=\delta(x_1)\big|Z_0=z\big)=c_z.\]
By Lemma \ref{lem:alpha}, $\alpha_{3,y_3}(z)$ is uniformly bounded below, while $\beta_3\gamma_3\cdots\beta_j\gamma_j$ depends only on the directed path (\ref{eqn:path}) and not on $z$, and thus is also uniformly bounded below since there are only finitely many $A_i$, and hence it suffices to choose finitely many such paths.

Similarly, $s_z$ also depends only on the directed path (\ref{eqn:path}), and thus is uniformly bounded above. In particular, it is possible to pick $s_z$ and $c_z$ such that $s=\sup_zs_z<\infty$ and $c=\inf_zc_z>0$.

Let $\tau$ be the first entry time into the state $\delta(x_1)$. By the above bound, we have $\P\big(\tau>s\big|Z_0=z\big)\le1-c$ for any initial state $z$ reachable from $\delta(x_1)$. Letting $\tau'$ and $Z'$ be independent copies of $\tau$ and $Z$, the time-homogeneous Markov property gives
\begin{align}
\frac{\P\big(\tau>(k+1)s\big|Z_0=z\big)}{\P\big(\tau>ks\big|Z_0=z\big)}
&=\P\big(\tau>(k+1)s\big|\tau>ks,Z_0=z\big)\nonumber\\
&=\sum_{z'}\P\big(\tau>(k+1)s\big|\tau>ks,Z_{ks}=z',Z_0=z\big)\nonumber\\[-.6cm]
&\hspace{2cm}\times\P\big(Z_{ks}=z'\big|\tau>ks,Z_0=z\big)\nonumber\\[.2cm]
&\le\sup_{z'}\P\big(\tau>(k+1)s\big|\tau>ks,Z_{ks}=z'\big)\nonumber\\
&=\sup_{z'}\P\big(\tau'>s\big|Z'_0=z'\big)\nonumber\\
&\le1-c.
\end{align}
By induction, $\P\big(\tau>ks\big|Z_0=z\big)\le(1-c)^k$ for all $k$. Dropping the condition on the initial distribution $Z_0$ for convenience, we have
\begin{align}
\E[\tau]
&=\sum_{k\in\Z^+}\P(\tau>k)=\sum_{k\in\Z^+}\sum_{0\le t<s}\P(\tau>ks+t)\nonumber\\
&\le\sum_{k\in\Z^+}\sum_{0\le t<s}\P(\tau>ks)\le\sum_{k\in\Z^+}s(1-c)^k=\frac sc<\infty.
\end{align}
In particular, $\E\big[\tau\big|Z=\delta(x_1)\big]<\infty$, so $\delta(x_1)$ is a positive recurrent state.
\end{proof}

\begin{lemma}
\label{lem:linearity}
The transition function $F$, considered as an operator on the real Banach space $\cM(\cP(S))$ of signed Radon measures on $\cP(S)$ with the total variation norm, is linear with operator norm $||F||=1$.
\end{lemma}
\begin{proof}
Linearity follows immediately from the fact that $F$ is defined as a finite sum of integrals.

For each $\mu\in\cM(\cP(S))$, let $\mu=\mu^+-\mu^-$ be the Hahn decomposition, so that $||\mu||=||\mu^+||+||\mu^-||=\mu^+(\cP(S))+\mu^-(\cP(S))$ by definition of the total variation norm.

If $\mu^+(\cP(S))=0$, then $F(\mu^+)(\cP(S))=0$ by linearity of $F$. Otherwise, let $c=\mu^+(\cP(S))$, so that $\frac1c\mu^+\in\cP(\cP(S))$. Since $F$ maps probability measures to probability measures, $||F(\mu^+)||=c||F\big(\frac1c\mu^+\big)||=c=||\mu^+||$, and similarly for $\mu^-$, hence by linearity of $F$ and the triangle inequality,
\begin{equation}
||F(\mu)||\le||F(\mu^+)||+||F(\mu^-)||=||\mu^+||+||\mu^-||=||\mu||.
\end{equation}
This shows that $||F||\le1$. Picking $\mu$ to be any probability measure gives $||F(\mu)||=1=||\mu||$, hence $||F||=1$.
\end{proof}

\begin{theorem}
\label{thm:convergence2}
In a positive, anchored hidden Markov model, irreducibility and aperiodicity of the restricted information chain $F|_R$ is a sufficient condition for convergence in distribution of the information state $Z_t$ to some discrete invariant measure $\mu_\infty\in\cP(R)\subset\cP(\cP(S))$.
\end{theorem}
\begin{proof}
In this case, the restricted chain $F|_R$ is irreducible, aperiodic and positive recurrent---that is, ergodic---and hence has a unique invariant probability measure $\mu_\infty$, which can be considered as an element of $\cP(\cP(S))$ supported on $R\subset\cP(S)$.

We will show that the information distribution $\mu_t$ converges in total variation norm to $\mu_\infty$. Fix $\mu_0\in\cP(\cP(S))$ and $\epsilon>0$, and pick $s$ such that $\lambda^s<\frac\epsilon3$, with $\lambda$ as defined in Proposition \ref{prop:discreteness}.

Let $\mu_s|_R$ be the restriction of the probability measure $\mu_s$ to $R$, which is a positive measure with total mass $m=\mu_s(R)\ge1-\lambda^s>0$, hence we can divide to obtain the probability measure $\mu'=\frac1m\mu_s|_R$ supported on $R$.

Since $F|_R$ is ergodic and $\mu'$ is a probability measure, $F^t(\mu')=(F|_R)^t(\mu')$ converges to $\mu_\infty$ in total variation norm, so pick $K$ such that for all $k\ge K$, $\bb F^k(\mu')-\mu_\infty\bb<\frac\epsilon3$.

For any $t\ge K+s$, we have the triangle inequality bound
\begin{align}
\bb\mu_t-\mu_\infty\bb
\le\bb\mu_t-F^{t-s}(m\mu')\bb
+\bb F^{t-s}(m\mu')-m\mu_\infty\bb
+\bb m\mu_\infty-\mu_\infty\bb.\label{eqn:triangleinequality}
\end{align}
By Lemma \ref{lem:linearity}, $F$ is linear with operator norm 1, so by Proposition \ref{prop:discreteness},
\begin{align}
\bb\mu_t-F^{t-s}(m\mu')\bb
&\le||F||^{t-s}\bb\mu_s-\mu_s|_R\bb\nonumber\\
&=\bb \mu_s|_{R^c}\bb=\mu_s(R^c)\le\lambda^s<\tfrac\epsilon3.
\end{align}
Since $t-s\ge K$, again using linearity and the fact that $m\le 1$,
\begin{equation}
\bb F^{t-s}(m\mu')-m\mu_\infty\bb=m\bb F^{t-s}(\mu')-\mu_\infty\bb<\tfrac\epsilon3.
\end{equation}
Finally, again using Proposition \ref{prop:discreteness}
\begin{equation}
\bb m\mu_\infty-\mu_\infty\bb=(1-m)\bb \mu_\infty\bb=1-\mu_s(R)\le\lambda^s<\tfrac\epsilon3.
\end{equation}
We see that for all $t\ge K+s$, $\bb\mu_t-\mu_\infty\bb<\epsilon$, so the information state  $\mu_t$ converges to $\mu_\infty$ in total variation norm and therefore in distribution.
\end{proof}

\begin{conjecture}
The conditions of Theorem \ref{thm:convergence2} can be weakened to the case when the underlying chain is only ergodic.
\end{conjecture}
\begin{proof}[Idea of Proof]
Given an ergodic finite transition matrix $T$, some power $T^k$ will be positive, and thus we should still have convergence by taking time in $k$-step increments, since the information state will not fluctuate too much within those $k$ steps. The difficulty lies in the fact that the information state taken in $k$-step increments is not the same as the information state of the $k$-step chain.
\end{proof}

It is our belief that the information state converges in all but a small number of pathological examples; however, we are only able to prove it in the above cases. If the information state does not converge, then it does not make sense to consider a limiting expected entropy. However, it is possible that a Ces\`aro sum of the expected entropy converges, and the limsup and liminf will certainly exist. Alternatively, we could simply work with a finite time horizon.

\newpage
\subsection{Special Case\label{sec:specialcase}}

Continuing further within the general case has proven to be quite difficult, so we will restrict the remainder of our results to a special case, where there are two states and two observation processes with two possible observations each, and each observation process observes a different underlying state perfectly. Formally, $|S|=|V|=|O|=2$, and the transition and observation matrices are
\begin{equation}
T=\bigtwomatrix{a}{1-a}{1-b}{b}\quad M^{(0)}=\bigtwomatrix{1}{0}{p}{1-p}\quad M^{(1)}=\bigtwomatrix{1-q}{q}{0}{1}.
\end{equation}
In order to exclude trivial cases, we will require that the parameters $a$, $b$, $p$ and $q$ are contained in the open interval $(0,1)$, and $a+b\ne1$.

Note that this special case exactly corresponds to the problem of searching for a moving target studied by MacPhee and Jordan \cite{macphee}, although our cost function, limiting expected entropy, is very different from theirs, expected cumulative sum of prescribed costs until the first zero-entropy observation.

We will begin by proving that given this restriction, the information state always converges in distribution, except for one case which is pathological in a sense that will be explained later. This proof is both a correction and an improvement of the proof given in \cite{honoursthesis}.

In the special case, the space $\cP(S)$ is a 1-simplex embedded in $\R^2$, which we can identify with the interval $[0,1]$, via equating each point $z\in[0,1]$ with the point $z\delta(0)+(1-z)\delta(1)\in\cP(S)$, so that $z$ represents the mass at 0 in the distribution.

By substituting these parameters into into Definition \ref{def:transitionfunction}, the transition function in the special case is
\begin{align}
F(\mu)&=\int_{A_0}\alpha_0(z)\delta\big(r_0(z)\big)d\mu(z) +\int_{A_1}\alpha_1(z)\delta\big(r_1(z)\big)d\mu(z)\label{eqn:specialcase}\\
&\hspace{1cm}+\int_{A_0}(1-\alpha_0(z))d\mu(z)\delta(0)+\int_{A_1}(1-\alpha_1(z))d\mu(z)\delta(1),\nonumber
\end{align}
where:
\vspace{-\parskip}
\begin{itemize}
\item $\alpha_0(z)=(1-p)(a+b-1)z+1-b+pb$;
\item $\alpha_1(z)=(1-q)(1-a-b)+b+q-bq$;
\item $r_0(z)=\displaystyle\frac{(a+b-1)z+1-b}{\alpha_0(z)}$; \rm{and}
\item $r_1(z)=\displaystyle\frac{q(a+b-1)z+q-qb}{\alpha_1(z)}$.
\end{itemize}
Note that the second line of (\ref{eqn:specialcase}) consists of two point masses at 0 and 1, which is a feature of the anchoredness condition. In the special case, it allows us to represent the location of masses by only two $r$-functions.

We will continue to use the symbols $\alpha_0$, $\alpha_1$, $r_0$ and $r_1$ in their meaning above for the remainder of our discourse. Note that we have simplified the notation, in that $r_0$ represents the $r$-function $r_{0,0}$, while $r_{0,1}$ does not appear since $r_{0,1}(z)=0$ identically, and similarly for the other symbols.

Since the special case satisfies the conditions of positivity and anchoredness, Proposition \ref{prop:positiverecurrent} applies, so irreducibility and aperiodicity of the information chain are sufficient for ergodicity. We now show that this occurs in all but two exceptional cases:
\vspace{-\parskip}
\begin{itemize}
\item Case 1: Each orbit is contained entirely in the policy region which maps to that orbit, that is, $R_0\subseteq A_0$ and $R_1\subseteq A_1$. Note that by (\ref{eqn:specialcase}), the $\alpha$-functions are always strictly positive, so we must have $R_0=\big\{0,r(0),r^2(0),\ldots\big\}$ and $R_1=\big\{1,r(1),r^2(1),\ldots\big\}$.
\item Case 2: The orbits alternate periodically between policy regions, that is, $\big\{1,r(0),r^2(1),r^3(0),\ldots\big\}\in A_0$ and $\big\{0,r(1),r^2(0),r^3(1),\ldots\big\}\in A_1$, where $r$ is the combined $r$-function $r(z)=r_0(z)\I_{A_0}(z)+r_1(z)\I_{A_1}(z)$.
\end{itemize}
\vspace{-\parskip}
Let Case 0 denote the general case when neither Case 1 nor Case 2 occurs.

\begin{lemma}
\label{lem:irreducible}
The chain $F|_R$ has only one irreducible recurrent class, except in Case 1, where it splits into two irreducible recurrent classes, both of which are positive recurrent.
\end{lemma}
\begin{proof}
By Proposition \ref{prop:positiverecurrent}, without loss of generality, assume that the state 0 is positive recurrent. 

If 1 is reachable from 0, that is, there is some $t$ such that $(F|_R)^t(0,1)>0$, then 0 is also reachable from 1 since 0 is recurrent, hence 0 and 1 are in the same irreducible recurrent class. By (\ref{eqn:specialcase}), either 0 or 1 is reachable from every state, so there cannot be any other recurrent classes.

Otherwise, if 0 is reachable from 1 but 1 is not reachable from 0, then 1 is transient, and furthermore, all of $R_1$ is transient since any $r^k(1)\in R_1$ is reachable only via the transient state 1, while all of $R_0$ is reachable from the recurrent state 0 and hence forms the only irreducible recurrent class.

Finally, if 0 and 1 are both unreachable from each other, then it must be the case that $R_0\subseteq A_0$ and $R_1\subseteq A_1$, in which case the chain splits into two irreducible classes $R_0$ and $R_1$, both of which are positive recurrent by the argument in Proposition \ref{prop:positiverecurrent}.
\end{proof}

\begin{lemma}
\label{lem:aperiodic}
The recurrent classes are aperiodic, except in Case 2, where the chain consists of a single irreducible recurrent class with period 2.
\end{lemma}
\begin{proof}
If any recurrent class is periodic, then at least one of 0 or 1 is recurrent and periodic; without loss of generality suppose it is 0. Since 0 is periodic, it cannot reach itself in 1 step, so must be contained in $A_1$, thus 1 is reachable from 0 and hence is in the same irreducible recurrent class. Note that by the same argument, $1\in A_0$, thus 0 reaches itself in 2 steps and hence the period must be 2.

This means 0 cannot reach itself in an odd number of steps, so $r^k(0)\in A_1$ when $k$ is even and $r^k(0)\in A_0$ when $k$ is odd, and similarly for the orbit of 1, which is the only possibility of periodicity.
\end{proof}

Thus, there is one exception in which the information chain is periodic, and another in which it is reducible. As will be evident, both exceptions can be remedied. We begin by giving a class of policies under which reducibility cannot occur. This class of policies is simple enough to be easily analysed, and as will be conjectured later, always includes the optimal policy in any hidden Markov model within the special case.

\begin{definition}
\label{def:thresholdpolicy}
A policy $g$ is called a \textbf{Threshold Policy} if its preimages $A_0=g^{-1}\{0\}$ and $A_1=g^{-1}\{1\}$ are both intervals.
\end{definition}

A threshold policy is indeed given by a threshold, since there must be some boundary point between $A_0$ and $A_1$, such that one observation process is used on one side and the other is used on the other side. Outside the threshold case, it is unclear what the equivalent definition would be, since the concept of an interval does not generalise easily to higher dimensions.

\begin{lemma}
\label{lem:r0bigger}
The linear fractional transformations $r_0$ and $r_1$ satisfy the inequality $r_0(z)>r_1(z)$ for all $z\in[0,1]$.
\end{lemma}
\begin{proof}
We can write $r_0(z)=\frac{(a+b+1)z+1-b}{p((1-a)z+b(1-z))+(a+b-1)z+1-b}$. Note that the coefficient $(1-a)z+b(1-z)$ of $p$ is strictly positive and in the denominator, while $r_1(z)$ is exactly the same with $q^{-1}$ instead of $p$. Since $p,q\in(0,1)$, $p<1<q^{-1}$, hence $r_0(z)>r_1(z)$.
\end{proof}

\begin{lemma}
\label{lem:monotonic}
The linear fractional transformations $r_0$ and $r_1$ are both strictly increasing when $a+b>1$ and both strictly decreasing when $a+b<1$.
\end{lemma}
\begin{proof}
The derivative of $r_0$ is $r_0'(z)=\frac{p(a+b-1)}{((1-p)(a+b-1)z+1-b+pb)^2}$, which is positive everywhere if $a+b>1$ and negative everywhere if $a+b<1$. The same holds for $r_1$, since it is identical with $q^{-1}$ instead of $p$.
\end{proof}

\begin{lemma}
\label{lem:fixedpoint}
The linear fractional transformations $r_0$ and $r_1$ have unique fixed points $\eta_0$ and $\eta_1$, which are global attractors of their respective dynamical systems.
\end{lemma}
\begin{proof}
Split the interval $[0,1]$ into open subintervals at its interior fixed points, noting that by inspection, the two boundary points are not fixed. Since linear fractional transformations have at most two fixed points, there must be between one and three such subintervals, all non-degenerate.

By continuity, in any such subinterval $I$, either $r_0(z)>z$ for all $z\in I$, or $r_0(z)<z$ for all $z\in I$. Since $r_0(0)>0$, $r_0(z)>z$ for all points $z$ in the leftmost subinterval, and since $r_0(1)<1$, $r_0(z)<z$ for all points $z$ in the rightmost subinterval. This means that there are at least two subintervals.

If there were three subintervals, $I_1$, $I_2$ and $I_3$ in that order, then either $r_0(z)>z$ for all $z\in I_2$, or $r_0(z)<z$ for all $z\in I_2$. In the first case, the fixed point between $I_1$ and $I_2$ is attracting on the left and repelling on the right, and in the second case, the fixed point between $I_2$ and $I_3$ is repelling on the left and attracting on the right. However, such fixed points can only occur for parabolic linear fractional transformations with only one fixed point. This is a contradiction, hence there cannot be three subintervals.

Thus, there are exactly two subintervals, and exactly one fixed point, which must then be an attractor across all of $[0,1]$.
\end{proof}

\begin{lemma}
\label{lem:eta0bigger}
The fixed points satisfy $\eta_1<\eta_0$.
\end{lemma}

\begin{proof}
By Lemma \ref{lem:r0bigger}, $r_1(\eta_0)<r_0(\eta_0)=\eta_0$.

First consider the case $a+b>1$. Applying Lemma \ref{lem:monotonic} $k$ times gives $r_1^{k+1}(\eta_0)<r_1^k(\eta_0)$, so the orbit of $\eta_0$ under $r_1$ is monotonically decreasing, but it also converges to $\eta_1$ by Lemma \ref{lem:fixedpoint}, hence $\eta_1<\eta_0$.

In the remaining case $a+b<1$, suppose $\eta_1\ge\eta_0$. Then by Lemma \ref{lem:monotonic}, $\eta_1=r_1(\eta_1)\le r_1(\eta_0)<\eta_0$, which is a contradiction, hence $\eta_1<\eta_0$.
\end{proof}

\begin{proposition}
\label{prop:exception}
The first exception to Theorem \ref{thm:convergence2}, Case 1, cannot occur under a threshold policy.
\end{proposition}
\begin{proof}
Suppose $R_0\subseteq A_0$, $R_1\subseteq A_1$, and the policy is threshold. Since $0\in R_0$ and $1\in R_1$, it follows that every point of $R_0$ is less than every point of $R_1$. Since $R_0$ and $R_1$ are the orbits of 0 and 1 under $r_0$ and $r_1$ respectively, by Lemma \ref{lem:fixedpoint}, they have limit points $\eta_0$ and $\eta_1$ respectively, hence $\eta_0\le\eta_1$. This contradicts Lemma \ref{lem:eta0bigger}, hence this cannot occur.
\end{proof}

The remaining exception is when the information chain is periodic with period 2, in which case the expected entropy oscillates between two limit points. The limiting expected entropy can still be defined in a sensible way, by taking the average, minimum or maximum of the two limit points, depending on which is most appropriate for the situation. Thus, for threshold policies, it is possible to define optimality without exception.

We conclude this section by writing down a closed form general expression for the limiting expected entropy.

\begin{proposition}
\label{prop:entropyformula}
Under the conditions of Theorem \ref{thm:convergence2}, that is, in Case 0, the limiting expected entropy of a policy is given by
\begin{equation}
H_\infty=\frac{C^{(0)}(H)+C^{(1)}(H)}{C^{(0)}(\mathbf1)+C^{(1)}(\mathbf1)},
\end{equation}
where, for $i\in\{0,1\}$:
\vspace{-\parskip}
\begin{itemize}
\item $H(z)=-z\log z-(1-z)\log(1-z)$ is the entropy function and $\mathbf1(z)=1$ is the constant function with value 1;
\item $r(z)=r_0(z)\I_{A_0}(z)+r_1(z)\I_{A_1}(z)$ and $\alpha(z)=\alpha_0(z)\I_{A_0}(z)+\alpha_1(z)\I_{A_1}(z)$ are the combined $r$-function and combined $\alpha$-function respectively, with $r_0$, $r_1$, $\alpha_0$ and $\alpha_1$ defined as in (\ref{eqn:specialcase});
\item $z^{(i)}_k\in\cP(S)=[0,1]$ are defined by the recursion
\begin{equation}
z^{(i)}_0=i,\qquad z^{(i)}_{k+1}=r(z_k^{(i)});
\end{equation}
\item $c^{(i)}_k\in\R$ are defined by the recursion
\begin{equation}
c^{(i)}_0=1,\qquad c^{(i)}_{k+1}=\alpha(z_k^{(i)})c_k^{(i)};
\end{equation}
\item $C^{(i)}:\cC(\cP(S))\rightarrow\R$ is the linear functional
\begin{equation}
C^{(i)}(f)=\sum_{k\in\Z^+}c_k^{(1-i)}(1-\alpha_i(z_k^{(1-i)}))\I_{\big\{z_k^{(1-i)}\in A_i\big\}}\sum_{k\in\Z^+}c_k^{(i)}f(z_k^{(i)}).
\end{equation}
\end{itemize}
\end{proposition}

\begin{proof}
By Theorem \ref{thm:convergence2}, a unique invariant probability measure $Z_\infty$ exists, so it suffices to show that
\begin{equation}
Z_\infty=\frac{C^{(0)}(\delta)+C^{(1)}(\delta)}{C^{(0)}(\mathbf1)+C^{(1)}(\mathbf1)},
\end{equation}
where $\delta:\cP(S)\rightarrow\cM(\cP(S))$ is the Dirac measure, and $C^{(i)}$ is extended in the obvious way to a linear functional $\cL(\cP(S),\cM(\cP(S)))\rightarrow\cM(\cP(S))$.

By Proposition \ref{prop:discreteness}, the invariant measure is supported in the combined orbit set $R$, and $z^{(i)}_k$ is the only point which can make a one-step transition under the restricted information chain $F|_R$ to $z^{(i)}_{k+1}$, with probability $\alpha(z^{(i)}_k)$. Since the masses at $z^{(i)}_k$ and $z^{(i)}_{k+1}$ are invariant, we must have
\begin{equation}
\P\big(Z_\infty=z^{(i)}_{k+1}\big)=\alpha(z^{(i)}_k)\P\big(Z_\infty=z^{(i)}_k\big).
\end{equation}
It then follows that for some constants $B^{(0)},B^{(1)}\in\R$,
\begin{equation}
Z_\infty=B^{(0)}\sum_{k\in\Z^+}c^{(0)}_k\delta(z^{(0)}_k) +B^{(1)}\sum_{k\in\Z^+}c^{(1)}_k\delta(z^{(1)}_k).
\end{equation}
The mass at $z^{(i)}_k$ is $B^{(i)}c^{(i)}_k$, which makes a transition under $F|_R$ to 0 in one step with probability $(1-\alpha_0(z^{(i)}_k))\I_{\{z^{(i)}_k\in A_0\}}$. Since $Z_\infty$ is an invariant measure, mass is conserved at $z_0^{(0)}=0$, hence
\begin{equation}
\label{eqn:conservationofmass}
\sum_{i\in\{0,1\}}B^{(i)}\sum_{k\in\Z^+}c^{(i)}_k(1-\alpha(z^{(i)}_k))\I_{\big\{z^{(i)}_k\in A_0\big\}}=B^{(0)}.
\end{equation}
Since $c^{(i)}_k\rightarrow0$ as $k\rightarrow\infty$, by telescoping series,
\begin{equation}
\label{eqn:telescope}
\sum_{k\in\Z^+}c^{(i)}_k(1-\alpha(z^{(i)}_k))=\sum_{k\in\Z^+}(c^{(i)}_k-c^{(i)}_{k+1})=1.
\end{equation}
Multiplying the right hand side of (\ref{eqn:conservationofmass}) by (\ref{eqn:telescope}) with $i=0$, then collecting coefficients of $B^{(0)}$ and $B^{(1)}$, yields
\begin{equation}
\label{eqn:conservationofmass2}
B^{(0)}\sum_{k\in\Z^+}c^{(0)}_k(1-\alpha(z^{(0)}_k))\I_{\big\{z^{(0)}_k\in A_1\big\}}=B^{(1)}\sum_{k\in\Z^+}c^{(1)}_k(1-\alpha(z^{(1)}_k))\I_{\big\{z^{(1)}_k\in A_0\big\}}.
\end{equation}
Note that conservation of mass at $\delta(0)=1$ is now automatic, since mass is conserved on all of $R$ and at every other point in $R$. The second equation comes from requiring $Z_\infty$ to have total mass 1, that is,
\begin{equation}
\label{eqn:totalmassone}
B^{(0)}\sum_{k\in\Z^+}c^{(0)}_k+B^{(1)}\sum_{k\in\Z^+}c^{(1)}_k=1.
\end{equation}
The solution to (\ref{eqn:conservationofmass2}) and (\ref{eqn:totalmassone}) is exactly the required result. We note that the denominator is zero only when $z^{(0)}_k\in A_0$ and $z^{(1)}_k\in A_1$ for all $k\in\Z^+$, which is exactly the excluded case $R_0\subseteq A_0$ and $R_1\subseteq A_1$.
\end{proof}

This proof can be generalised easily to the case of more than two observation processes, as long as each one is anchored with only two states.

\newpage
\section{Computational Results}
\label{chap:computation}
\subsection{Limiting Entropy}
We now present computational results, again in the special case, which will illustrate the nature of the optimal policy in the case of minimising limiting expected entropy. Since any such discussion must cover how best to calculate the limiting expected entropy given a hidden Markov model and a policy, this is a natural place to start.

The simpliest approach is to use the formula in Proposition \ref{prop:entropyformula}, noting that each of $C^{(0)}(H)$, $C^{(1)}(H)$, $C^{(0)}(\mathbf 1)$ and $C^{(1)}(\mathbf 1)$ is a product of two infinite series, each of which is bounded by a geometric series and hence well-approximated by truncation. Specifically, we write the limiting expected entropy as
\[H_\infty=\frac{C_1H_0+C_0H_1}{C_1I_0+C_0I_1},\]
where:
\vspace{-\parskip}
\begin{itemize}
\item $I_i=\displaystyle\sum_{k\in\Z^+}c_k^{(i)}$;
\item $H_i=\displaystyle\sum_{k\in\Z^+}c_k^{(i)}H(z_k^{(i)})$;
\item $C_i=\displaystyle\sum_{k\in\Z^+}c_k^{(i)}\big(1-\alpha_{1-i}(z_k^{(i)})\big)\I_{\big\{z_k^{(i)}\in A_{1-i}\big\}}$.
\end{itemize}
We can simplify the calculation slightly by recursively updating $c_k^{(i)}$ and $z_k^{(i)}$, storing them as real-valued variables $c_i$ and $z_i$.

\newpage
\begin{algorithm}
\label{alg:entropy}
Estimation of limiting expected entropy.
\vspace{-\parskip}
\begin{enumerate}
\item Define as functions the entropy $H(z)=-z\log z-(1-z)\log(1-z)$, the mixed $r$-function $r(z)=\I_{A_0}(z)r_0(z)+\I_{A_1}(z)r_1(z)$ and the mixed $\alpha$-function $\alpha(z)=\I_{A_0}(z)\alpha_0(z)+\I_{A_1}(z)\alpha_1(z)$;
\item Pick a large number $N$;
\item Initialise variables $\widehat C_0$, $\widehat C_1$, $\widehat H_0$, $\widehat H_1$, $\widehat I_0$ and $\widehat I_1$ to 0, $c_0$ and $c_1$ to 1, $z_0$ to 0 and $z_1$ to 1;
\item Repeat the following $N+1$ times:
\begin{enumerate}
\item Add $c_0$ to $\widehat I_0$ and $c_1$ to $\widehat I_1$;
\item Add $c_0H(z_0)$ to $\widehat H_0$ and $c_1H(z_1)$ to $\widehat H_1$;
\item If $z_0\in A_1$, then add $c_0(1-\alpha_1(z_0))$ to $\widehat C_0$, and if $z_1\in A_0$, then add $c_1(1-\alpha_0(z_1))$ to $\widehat C_1$;
\item Let $z_0=r(z_0)$ and $z_1=r(z_1)$;
\item Multiply $c_0$ by $\alpha(z_0)$ and $c_1$ by $\alpha(z_1)$;
\end{enumerate}
\item The limiting expected entropy $H_\infty$ is estimated by
\[\widehat H_\infty(N)=\frac{\widehat C_1\widehat H_0+\widehat C_0\widehat H_1}{\widehat C_1\widehat I_0+\widehat C_0\widehat I_1}.\]
\end{enumerate}
\end{algorithm}

\begin{proposition}
\label{prop:errorbound}
The estimate $\widehat H_\infty(N)$ satisfies the error bound
\[\big|\widehat H_\infty(N)-H_\infty\big|\le\frac{16\alpha^N}{(1-\alpha)^4Q(N)^2},\]
where $Q(N)=\widehat C_1\widehat I_0+\widehat C_0\widehat I_1$ is the denominator of $\widehat H_\infty(N)$, and
\[\alpha=\sup_{i,z}\alpha_i(z)=1-\min\Big\{b(1-p),(1-a)(1-p),(1-b)(1-q),a(1-q)\Big\}.\]
\end{proposition}
\begin{proof}
The formula for $\alpha$ follows from Definition \ref{def:transitionfunction}, since the $\alpha$-functions are linear and hence their maxima occur at the endpoints $z=0$ or $z=1$. Note that $\alpha<1$ follows from the requirement that $a,b,p,q\in(0,1)$. Since $Q(N)$ is monotonic increasing, this implies that the error bound is finite and vanishes as $N\rightarrow\infty$.

Since each series has non-negative summands, each truncated series is smaller than the untruncated series, hence
\begin{align*}
\bigg|H_\infty-\frac{C_1H_0+C_0H_1}{\widehat C_1\widehat I_0+\widehat C_0\widehat I_1}\bigg|
&=\frac{(C_1H_0+C_0H_1)(C_1I_0+C_0I_1-\widehat C_1\widehat I_0-\widehat C_0\widehat I_1)}{(C_1I_0+C_0I_1)(\widehat C_1\widehat I_0+\widehat C_0\widehat I_1)}\\
&\le\frac{C_1H_0+C_0H_1}{(\widehat C_1\widehat I_0+\widehat C_0\widehat I_1)^2}\big((C_1I_0-\widehat C_1\widehat I_0)+(C_0I_1-\widehat C_0\widehat I_1)\big).
\end{align*}
The $k$th summand in each series is bounded by $c^{(i)}_k\le\alpha^k$, hence
\begin{align*}
C_1I_0-\widehat C_1\widehat I_0
&=I_0(C_1-\widehat C_1)+C_1(I_0-\widehat I_0)-(C_1-\widehat C_1)(I_0-\widehat I_0)\\
&\le\frac1{1-\alpha}\frac{\alpha^N}{1-\alpha}+\frac1{1-\alpha}\frac{\alpha^N}{1-\alpha}=\frac{2\alpha^N}{(1-\alpha)^2}.
\end{align*}
The same bound holds for $C_0I_1-\widehat C_0\widehat I_1$, hence
\[\bigg|H_\infty-\frac{C_1H_0+C_0H_1}{\widehat C_1\widehat I_0+\widehat C_0\widehat I_1}\bigg|
\le\frac{8\alpha^N}{(1-\alpha)^4Q(N)^2}.\]
Similarly, since $Q(N)=\widehat C_1\widehat I_0+\widehat C_0\widehat I_1\le2/(1-\alpha)^2$,
\begin{align*}
\bigg|\widehat H_\infty(N)-\frac{C_1H_0+C_0H_1}{\widehat C_1\widehat I_0+\widehat C_0\widehat I_1}\bigg|
&=\frac{(C_1H_0-\widehat C_1\widehat H_0)+(C_0H_1-\widehat C_0\widehat H_1)}{\widehat C_1\widehat I_0+\widehat C_0\widehat I_1}\\
&\le\frac{4\alpha^N}{(1-\alpha)^2Q(N)}\le\frac{8\alpha^N}{(1-\alpha)^4Q(N)^2}.
\end{align*}
Combining via the triangle inequality gives the required bound.
\end{proof}

Note that this bound depends only on the quantities $\alpha$ and $Q(N)$, which are easily calculated. In particular, this allows us to prescribe a precision and calculate the limiting expected entropy to within that precision by running the algorithm with unbounded $N$ and terminating when error bound reaches the desired precision. Furthermore, since $Q(N)$ appears only in the denominator and grows monotonically with $N$, we can replace it with $Q(N_0)$ for some small, fixed $N_0$ to calculate a prior a sufficient number of steps for any given precision.

\begin{example}
\label{eg:errorbound}
In later simulations, we will pick $a,b,p,q\in[0.025,0.975]$, so that $\alpha\le0.999375$. This gives an error bound of
\[\big|\widehat H_\infty(N)-H_\infty\big|\le(16\times400^4)Q(N)^{-2}(0.999375)^N.\]
While the constant appears daunting at first glance, solving for a prescribed error of $10^{-10}$ gives
\[\log 16+4\log 400-2\log Q(N)+N\log 0.999375\le-10\log 10.\]
Hence, we require
\[N\ge79598+3199\,\big|\!\log Q(N)\big|.\]
For any realistic value of $Q(N)$, this is easily within computational bounds, as each iteration requires at most 36 arithmetic operations, 2 calls to the policy function, and 4 calls to the logarithm function.
\end{example}

An alternative approach to estimating limiting expected entropy would be to simulate observation of the hidden Markov model under the given policy. The major drawback of this method is that it requires working with the information state $Z_t$, which takes values in $\cP(\cP(S))=\cP([0,1])$, the set of probability measures on the unit interval, which is an infinite dimensional space.

One possibility is to divide $[0,1]$ into subintervals and treat each subinterval as a discrete atom, but this produces a very imprecise result. Even using an unrealistically small subinterval width of $10^{-6}$, the entropy function has a variation of over $10^{-5}$ across the first subinterval, restricting the error bound to this order of magnitude regardless of the number of iterations. In comparison, Example \ref{eg:errorbound} shows that the direct estimation method has greatly superior performance.

An improvement is to use the fact that the limiting distribution $Z_\infty$ is discrete, and store a list of pairs containing the locations and masses of discrete points. Since any starting point moves to either 0 or 1 in one step, at the $N$th iteration, the list of points must contain at least the first $N$ points in the orbit of either 0 or 1. Each such point requires a separate calculation at each iteration, and thus the number of computations is $O(N^2)$ rather than $O(N)$ as for Algorithm \ref{alg:entropy}.

Since the number of iterations $N$ corresponds to the last point in the orbit of 0 or 1 which is represented, for any given $N$, this method differs from the direct computation method only in the masses on these points, thus we would expect the relationship between precision and number of iterations to be similar. Since the simulation method has quadratically growing number of computations, this would suggest that it is slower than the direct computation method, and indeed, this is also indicated by empirical trials.

We will use the direct computation method of estimating limiting expected entropy for all of our results.

\newpage
\subsection{Optimal Threshold Policy}
The problem of finding the policy which minimises limiting expected entropy is made much easier by restricting the search space to the set of threshold policies, as these can be represented by a single number representing the threshold, and a sign representing which observation process is used on which side of the threshold.

The simplest approach is to pick a collection of test thresholds uniformly in $[0,1]$, either deterministically or randomly, and test the limiting expected entropy at these thresholds, picking the threshold with minimal entropy as the optimal threshold policy. However, this method is extremely inefficient. Proposition \ref{prop:discreteness} shows that the policy only matters on the countable set $R\subset[0,1]$, so moving the threshold does not change the policy as long as it does not move past a point in $R$.

\begin{figure}[H]
  \centering
  \includegraphics[width=\textwidth]{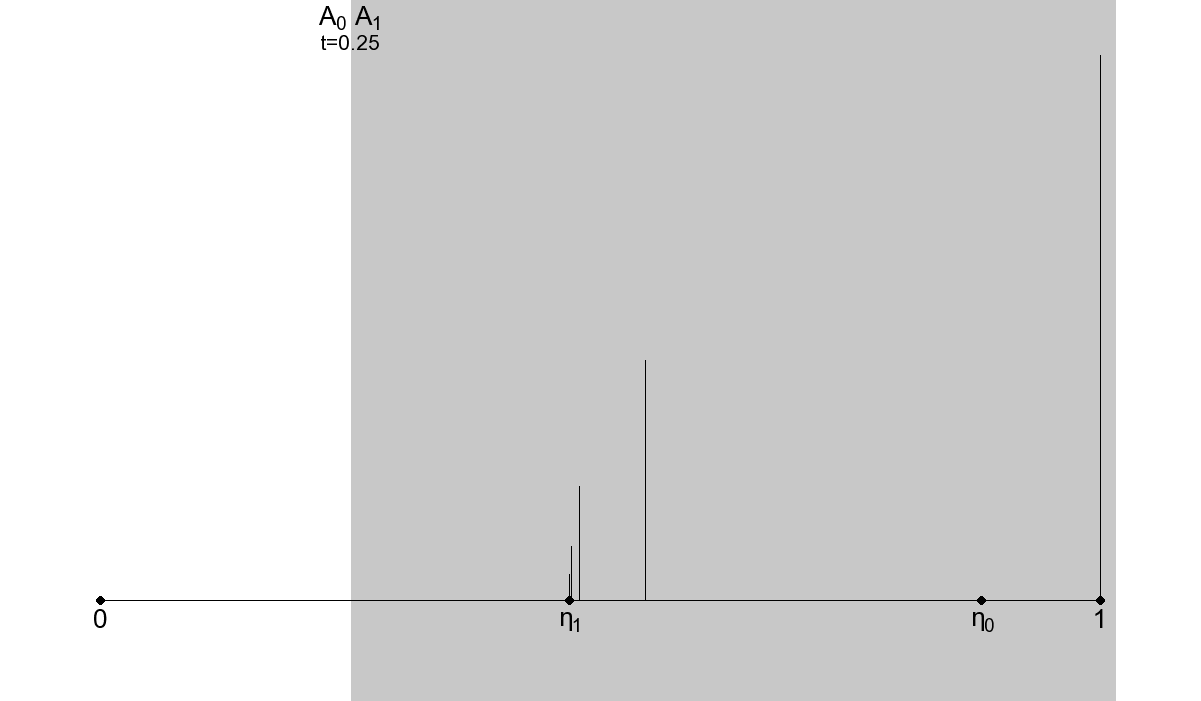}
  \large
  \caption[Example of a Region I policy]{Invariant measure for the unique Region I policy with $a=0.8$, $b=0.3$, $p=0.5$, $q=0.3$. Entropy is 0.3145. Under the evolution function, any mass eventually enters the grey region since it contains both accumulation points in its interior, after which it cannot escape, hence in the limit, there is zero mass in the orbit of 0, and the policy is equivalent the all-$A_1$ policy.}
  \label{fig:regionI}
\end{figure}

\begin{figure}[H]
  \centering
  \includegraphics[width=\textwidth]{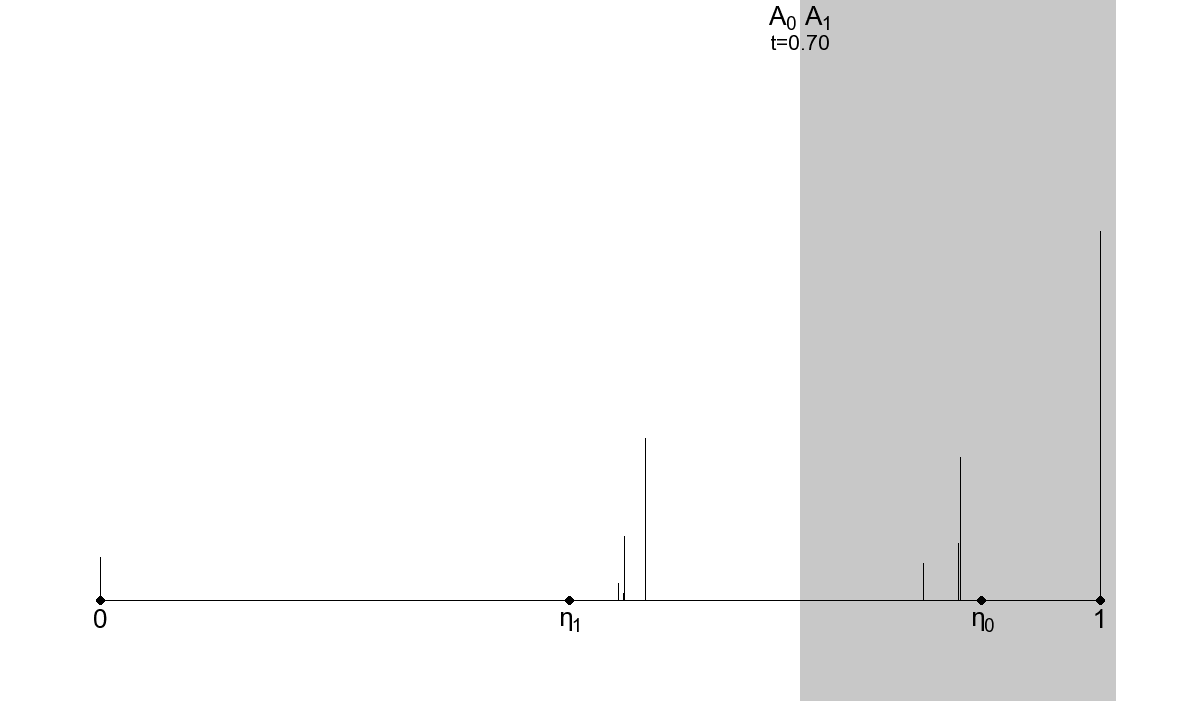}
  \large
  \caption[Example of a Region II policy]{Invariant measure for a typical Region II policy with $a=0.8$, $b=0.3$, $p=0.5$, $q=0.3$. Entropy is 0.3251. Note that the masses do not converge to the accumulation points.}
  \label{fig:regionII}
\end{figure}

\begin{figure}[H]
  \centering
  \includegraphics[width=\textwidth]{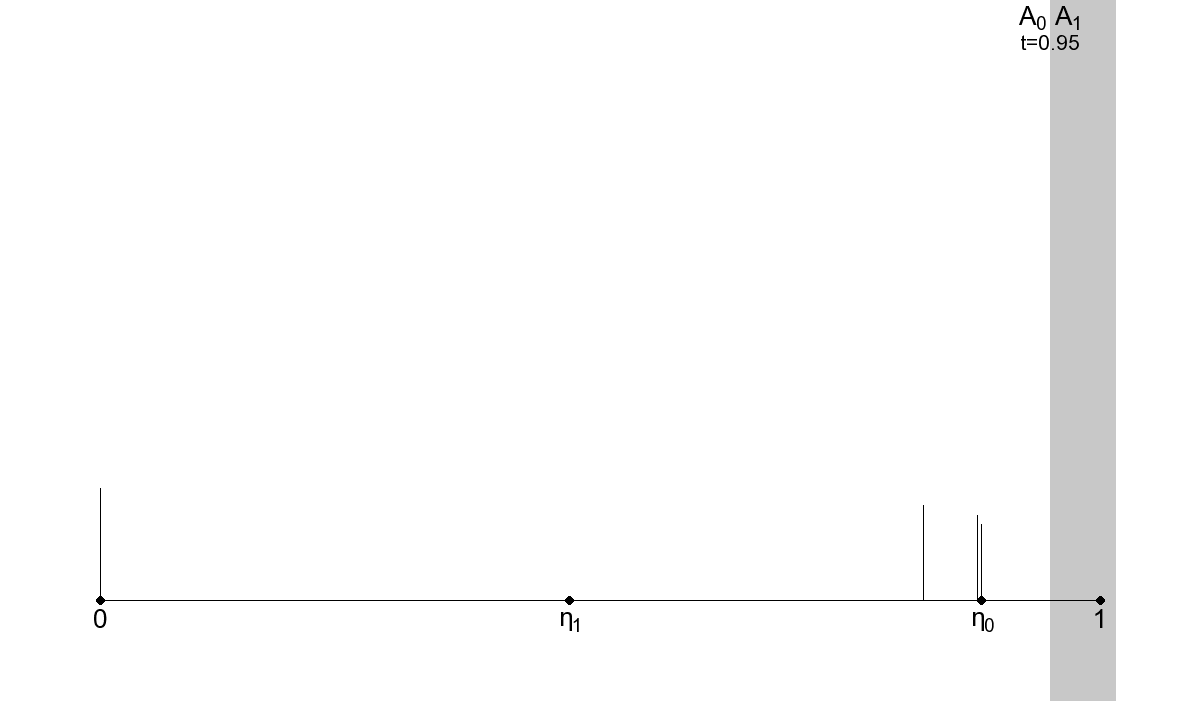}
  \large
  \caption[Example of a Region III policy]{Invariant measure for the unique Region III policy with $a=0.8$, $b=0.3$, $p=0.5$, $q=0.3$. Entropy is 0.3337. Under the evolution function, any mass eventually enters the white region since it contains both accumulation points in its interior, after which it cannot escape, hence in the limit, there is zero mass in the orbit of 1, and the policy is equivalent to the all-$A_0$ policy.}
  \label{fig:regionIII}
\end{figure}

\begin{figure}[H]
  \centering
  \includegraphics[width=\textwidth]{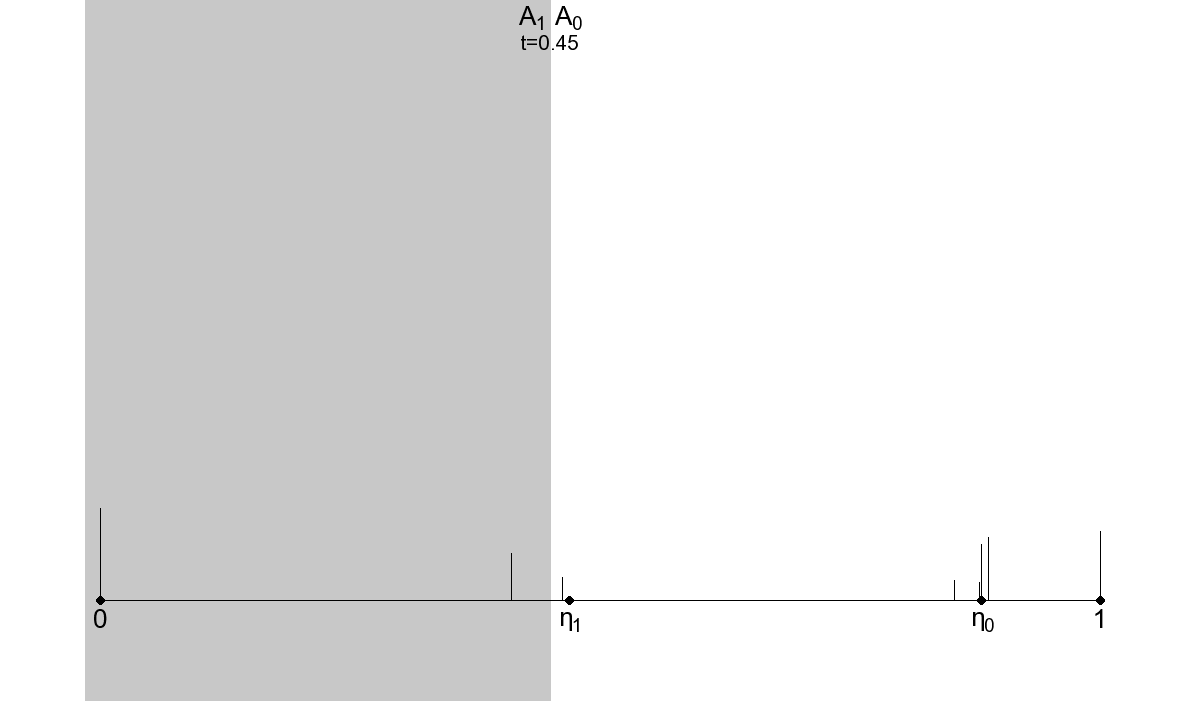}
  \large
  \caption[Example of a Region IV policy]{Invariant measure for a typical Region IV policy with $a=0.8$, $b=0.3$, $p=0.5$, $q=0.3$. Entropy is 0.3275. The orbit of 1 approaches $\eta_0$, while the orbit of 0 follows an approach sequence to $\eta_1$ for a finite number of steps (in this case 2 steps) before also approaching $\eta_0$.}
  \label{fig:regionIV}
\end{figure}

\begin{figure}[H]
  \centering
  \includegraphics[width=\textwidth]{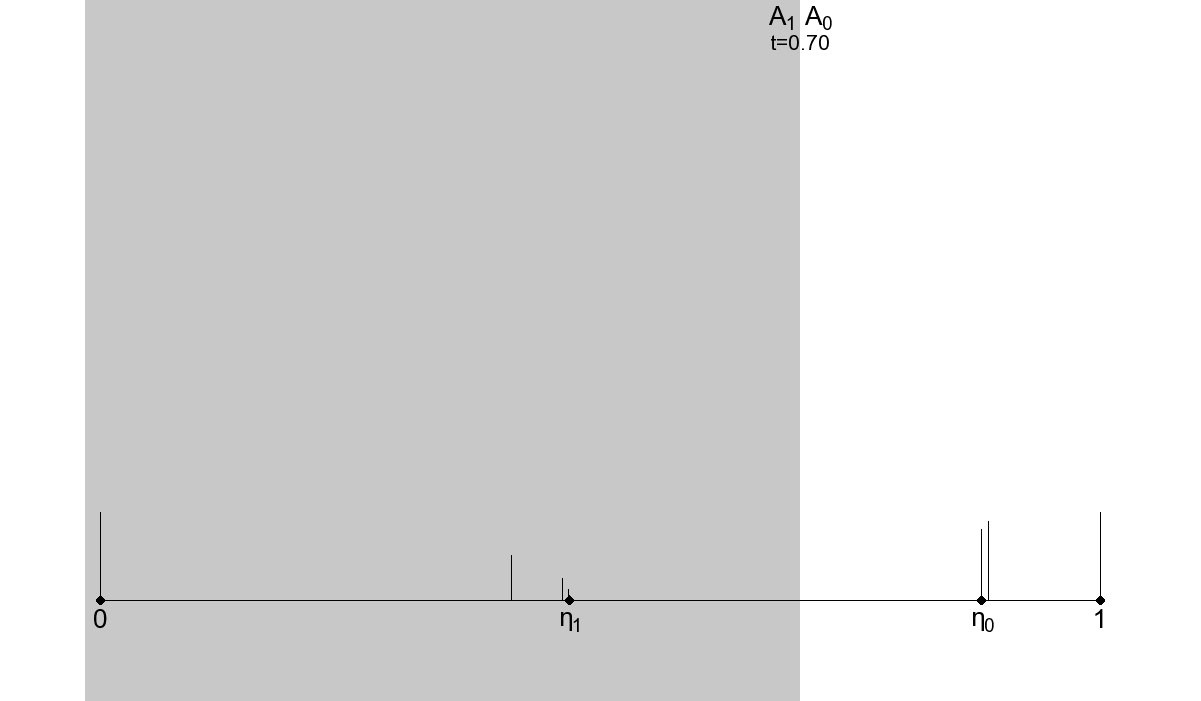}
  \large
  \caption[Example of a Region V policy]{Invariant measure for the unique Region V policy with $a=0.8$, $b=0.3$, $p=0.5$, $q=0.3$. Entropy is 0.3265. Note that the orbit of 0 converges to $\eta_1$ while the orbit of 1 converges to $\eta_0$.}
  \label{fig:regionV}
\end{figure}

\begin{figure}[H]
  \centering
  \includegraphics[width=\textwidth]{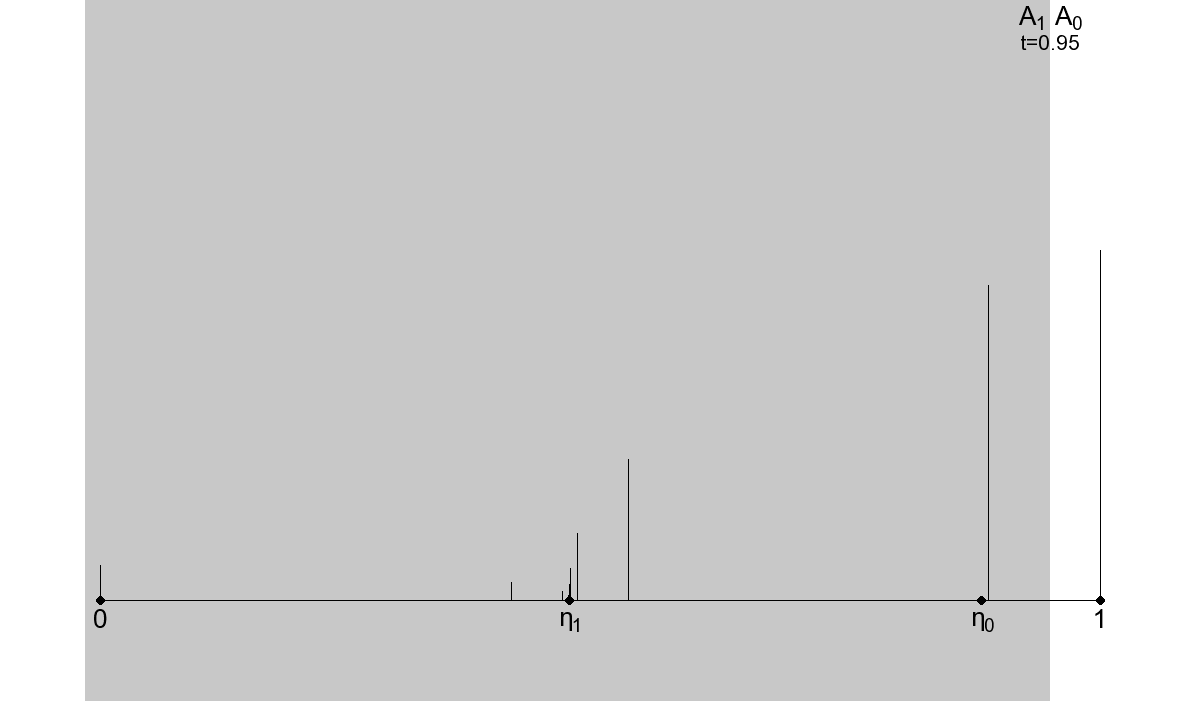}
  \large
  \caption[Example of a Region VI policy]{Invariant measure for a typical Region VI policy with $a=0.8$, $b=0.3$, $p=0.5$, $q=0.3$. Entropy is 0.3182. The orbit of 0 approaches $\eta_1$, while the orbit of 1 follows an approach sequence to $\eta_0$ for a finite number of steps (in this case 1 step) before also approaching $\eta_1$.}
  \label{fig:regionVI}
\end{figure}

As shown in Figures \ref{fig:regionI}--\ref{fig:regionVI}, points in $R$ tend to be quite far apart, and thus the naive approach will cause a large number of equivalent policies to be tested. On the other hand, points in $R$ close to the accumulation points are closely spaced, so even with a very fine uniform subset, some policies will be missed when the spacing between points in $R$ becomes less than the spacing between test points.

A better way is to decide on a direction in which to move the threshold, and select the test point exactly at the next point in the previous realisation of $R$ in the chosen direction, so that every threshold in between the previous point and the current point gives a policy equivalent to the previous policy. This ensures that each equivalence class of policies is tested exactly once, thus avoiding the problem with the naive method.

However, a new problem is introduced in that the set of test points depends on the iteration number $N$, which determines the number of points of $R$ that are considered. This creates a circular dependence, in that the choice of $N$ depends on the desired error bound, the error bound depends on the policy, and the set of policies to be tested depends on $N$. We can avoid this problem by adapting Proposition \ref{prop:errorbound} to a uniform error bound across all threshold policies.

\begin{proposition}
\label{prop:errorboundthreshold}
For a threshold policy and $N\ge L$, the error is
\[\big|\widehat H_\infty(N)-H_\infty\big|<\frac{16\alpha^N}{\bar\alpha^{2L}(1-\alpha)^6},\]
where $L$ is the smallest integer such that $r_0^L(0)>r_1^L(1)$, and
\[\bar\alpha=\inf_{i,z}\alpha_i(z)=1-\max\Big\{b(1-p),(1-a)(1-p),(1-b)(1-q),a(1-q)\Big\}.\]
\end{proposition}
\begin{proof}
First note that $L$ exists, since by Lemma \ref{lem:eta0bigger}, iterations of $r_0$ and $r_1$ converge to $\eta_0>\eta_1$ respectively.

Using Proposition \ref{prop:errorbound}, it suffices to prove that for $N\ge L$,
\[Q(N)=\widehat C_1\widehat I_0+\widehat C_0\widehat I_1\ge\bar\alpha^L(1-\alpha).\]
It is not possible for $\{0,r(0),\ldots,r^L(0)\}\subset A_0$ and $\{1,r(1),\ldots,r^L(1)\}\subset A_1$, since this would mean $r^L(0)=r_0^L(0)$ and $r^L(1)=r_1^L(1)$, which gives the ordering $0<r^L(1)<r^L(0)<1$, but $A_0$ and $A_1$ are intervals for a threshold policy.

Hence, either $z_\ell^{(0)}=r^\ell(0)\in A_1$ or $z_\ell^{(1)}=r^\ell(1)\in A_0$ for some $\ell\le L\le N$. If $z_\ell^{(0)}\in A_1$, then $\widehat C_0\ge c_\ell^{(0)}\big(1-\alpha_1(z_\ell^{(0)})\big)\ge\bar\alpha^\ell(1-\alpha)$. Since $\widehat I_1\ge c_0^{(1)}=1$, this gives $Q(N)\ge\bar\alpha^L(1-\alpha)$, as required. A similar argument holds in the case $z_\ell^{(1)}\in A_0$.
\end{proof}

The existence of this uniform bound for $Q(N)$ in the threshold case is closely related to Proposition \ref{prop:exception}, which states that the exception Case 1, where $R_0\subset A_0$ and $R_1\subset A_1$, cannot occur in a threshold policy. In this exceptional case, Proposition \ref{prop:entropyformula} does not hold, as the denominator is zero, and hence $Q(N)=0$ for all $N$. The fact that this cannot occur in a threshold policy is the key ingredient of this uniform bound.

Now that we have an error bound which does not depend on the policy, we can determine a uniform number of iterations $N$ that will suffice for estimating the limiting expected entropy for all threshold policies. This reduces the search space to a finite one, as each point in the orbits of 0 and 1 must be in one of the two policy regions, hence, there are at most $2^{2N+2}$ policies. Most of these will not be threshold policies, but since orbit points need not be ordered, there is no obvious bound on the number of threshold policies that need to be checked. Simulation results later in this section will show that in most cases, the number of such policies is small enough to be computationally feasible.

\begin{definition}
\label{def:orientation}
The \textbf{Orientation} of a threshold policy is the pair of whether $A_0$ is to the left or right of $A_1$, and whether the threshold is included in the left or right interval. Let $[A_0)[A_1]$, $[A_0](A_1]$, $[A_1)[A_0]$ and $[A_1](A_0]$ denote the four possibilities, with the square bracket indicating inclusion of the threshold and round bracket indicating exclusion.
\end{definition}

Our strategy for simplifying the space of threshold policies that need to be considered is to note that the policy only matters on $R$, the support of the invariant measure. Although $R$ depends on the policy, for a given orientation and threshold $t$, any policy with the same orientation and some other threshold $t'$ such that no points of $R$ lie between $t$ and $t'$ is an equivalent policy, in that sense that the invariant measure is the same, since no mass exists in the region where the two policies differ.

Thus, for each orientation, we can begin with $t=1^+$, and at each iteration, move the threshold left past the next point in $R$, since every threshold in between is equivalent to the previous threshold. Although $R$ changes at each step, this process must terminate in finite time since we already showed that there are only finitely many policies for any given $N$, and by testing equivalence classes of policies only once, it is likely to that far fewer steps are required than the $2^{2N+2}$ bound.

Furthermore, since $R$ is a discrete set, every threshold policy has an interval of equivalent threshold policies, so we can assume without loss of generality that the threshold is contained in the interval to the right, that is, only test the orientations $[A_0)[A_1]$ and $[A_1)[A_0]$.

\begin{algorithm}
\label{alg:optimalthreshold}
Finding the optimal threshold policy.
\vspace{-\parskip}
\begin{enumerate}
\item Find $L$, the smallest integer such that $r_0^L(0)>r_1^L(1)$, by repeated application of $r_0$ and $r_1$ to 0 and 1 respectively;
\item Prescribe an error $E$ and determine the number of iterations
\[N=\frac{\log E-\log 16+2L\log\bar\alpha+6\log(1-\alpha)}{\log\alpha};\]
\item Start with the policy $A_0=[0,t)$ and $A_1=[t,1]$ with $t=1^+$, that is, $A_0$ is the whole interval and $A_1$ is empty but considered as being to the right of $A_0$, and loop until $t=0$:
\begin{enumerate}
\item Run Algorithm \ref{alg:entropy}, and let the next value of $t$ be the greatest $z_k^{(i)}$ which is strictly less than the previous value of $t$;
\item If entropy is less than the minimum entropy for any policy encountered so far, record it as the new minimum;
\end{enumerate}
\item Repeat with $A_0$ considered to be on the right of $A_1$.
\end{enumerate}
\end{algorithm}

\begin{figure}[H]
  \centering
  \includegraphics[width=\textwidth]{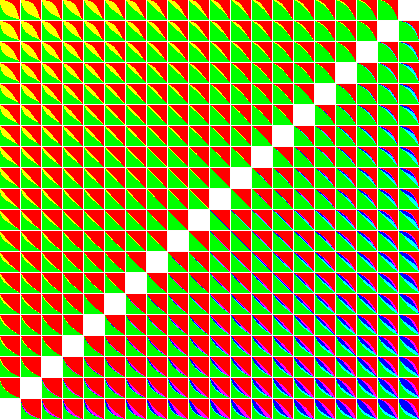}
  \large
  \caption[Location of the optimal threshold]{Location of the optimal threshold policy. The four-dimensional parameter space is represented with $a$ on the major horizontal axis increasing rightwards, $b$ on the major vertical axis increasing downwards, $p$ on the minor horizontal axis increasing rightwards and $q$ on the minor vertical axis increasing downwards. Each parameter is sampled by $\{0.025,0.075,\ldots,0.975\}$, omitting the trivial cases when $a+b=1$, resulting in 152000 sample points. Colours representing regions are defined in Figures \ref{fig:regionI}--\ref{fig:regionVI}.}
  \label{fig:thresholdregions}
\end{figure}

We now calculate the location of the optimal threshold for a range of values of the parameters $a$, $b$, $p$ and $q$. The results of this calculation, which will be presented in Figure \ref{fig:thresholdregions}, is the primary content of this section, as it will give an empirical description of the optimal threshold policy. Since we have not been able to prove optimality analytically, except in the symmetric case of Proposition \ref{prop:symmetric}, this empirical description will provide our most valuable insight into the problem.

In order to facilitate understanding, we will define six classes of threshold policies, depending on the location and orientation of the threshold. Diagrams of the invariant measure under each of these classes, presented in Figures \ref{fig:regionI}--\ref{fig:regionVI}, will demonstrate that these classes of thresholds are qualitatively different, which will further manifest itself in the next section.

We partition the interval into three subintervals with endpoints at the attractors $\eta_0$ and $\eta_1$, noting that this produces three regions which consist entirely of a single equivalence class of policy. We also assign colours to these regions to accommodate the data presented in Figure \ref{fig:thresholdregions}.

\begin{definition}
\label{def:thresholdregions}
The six \textbf{Threshold Regions} are defined by partitioning the set of all threshold policies first by whether the orientation is $[A_0)[A_1]$ or $[A_1)[A_0]$, then by the location of the threshold $t$ in relation to the accumulation points $\eta_1$ and $\eta_0$, with inclusion of endpoints defined more precisely overleaf. We number them using Roman numerals I through VI.
\end{definition}

Note that if either $t\le0$ or $t>1$, then either $A_0$ or $A_1$ occupies the whole interval, depending on the orientation. In particular, $[A_0)[A_1]$ with $t=0$ is equivalent to $[A_1)[A_0]$ with $t=1^+$, since in both cases $A_1$ is the whole interval, and similarly $[A_0)[A_1]$ with $t=1^+$ is equivalent to $[A_1)[A_0]$ with $t=0$.

Thus, the space of all possible threshold policies consists of two disjoint intervals $[0,1^+]$, each of whose endpoints is identified with an endpoint of the other interval, which is topologically equivalent to a circle. To be technically correct, we note that identifying 0 and $1^+$ does not present a problem, since we can simply extend the interval to $[0,1+\epsilon]$ for some small $\epsilon>0$, and identify 0 and $1+\epsilon$ instead. While this results in a subinterval of the circle corresponding to the same policy, this does not add additional complexity since every threshold policy has an interval of equivalent policies.

This topology is illustrated overleaf in Figure \ref{fig:thresholddomain}, followed by precise definitions and examples of the six threshold regions in Figures \ref{fig:regionI}--\ref{fig:regionVI}, and finally our computational results in Figure \ref{fig:thresholdregions}.

\vspace{3\parskip}
\begin{figure}[H]
  \centering
  \Large
  ~~~(I)~~~~~~~~~~~~~~(II)~~~~~~~~~~~~~(III)\\
  \includegraphics[width=\textwidth]{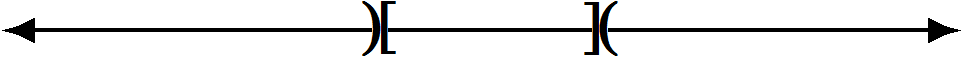}\\[-.5cm]
  0~~~~~~~~~~~~~~~~~~~$\eta_1$~~~~~~~~~~~~~~ $\eta_0$~~~~~~~~~~~~~~~~~1\\[.3cm]
  \includegraphics[width=\textwidth]{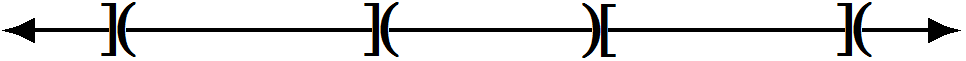}\\[-.5cm]
  (III)~~~~~~~~~~(IV)~~~~~~~~~~~~~~(V)~~~~~~~~~~~~~ (VI)~~~~~~~~~(I)~~\\[.8cm]
  \large
  \caption[Space of all threshold policies]{Space of all threshold policies. The top line represents the orientation $[A_0)[A_1]$, while the bottom line represents the orientation $[A_1)[A_0]$. The right end of the top line and the left end of the bottom line are both the policy $A_0=[0,1]$ and $A_1=\emptyset$, so we can paste them together, and similarly for the left end of the top line and the right end of the bottom line. Hence, we see that the set of threshold policies is topologically a circle.}
  \label{fig:thresholddomain}
\end{figure}
\vspace\parskip

\newpage
\textbf{Region I} (represented by red in Figure \ref{fig:thresholdregions}): $[A_0)[A_1]$ with $t<\eta_1$ or $[A_1)[A_0]$ with $t>1$. When $a+b>1$, every policy here is equivalent to the all-$A_1$ policy, as the mass in the orbit of 0 approaches 0.
\vspace{3\parskip}
\vspace\parskip

\newpage
\textbf{Region II} (yellow): $[A_0)[A_1]$ with $\eta_1\le t\le\eta_0$. This is the most difficult to understand of the threshold policies, as the orbits do not converge to the accumulation points $\eta_0$ and $\eta_1$, but rather, oscillate around the threshold $t$.
\vspace{3\parskip}
\vspace\parskip

\newpage
\textbf{Region III} (green): $[A_0)[A_1]$ with $t>\eta_0$ or $[A_1)[A_0]$ with $t\le0$. When $a+b>1$, every policy here is equivalent to the all-$A_0$ policy, since the mass in the orbit of 1 approaches 0.
\vspace{3\parskip}
\vspace\parskip

\newpage
\textbf{Region IV} (cyan): $[A_1)[A_0]$ with $0<t\le\eta_1$. For a policy in this region, the first finitely many points of the orbit of 0 belongs to $A_1$, while every other point lies in $A_0$. Note that $t\le0$ is included in Region III.
\vspace{3\parskip}
\vspace\parskip

\newpage
\textbf{Region V} (blue): $[A_1)[A_0]$ with $\eta_1<t<\eta_0$. When $a+b>1$, every policy here is equivalent to the policy $R_0\subset A_1$ and $R_1\subset A_0$.
\vspace{3\parskip}
\vspace\parskip

\newpage
\textbf{Region VI} (magenta): $[A_1)[A_0]$ with $\eta_0\le t\le1$. For a policy in this region, the first finitely many points of the orbit of 0 belongs to $A_1$, while every other point lies in $A_0$. Note that $t>1$ is included in Region I. The lack of symmetry with Region VI in terms of strict and non-strict inequalities is due to the choice that the threshold itself be included in the right region. 
\vspace{3\parskip}
\vspace\parskip

\vspace{3\parskip}
\vspace\parskip

The bound in Proposition \ref{prop:errorboundthreshold} is consistent with the empirical running time for the calculations used to generate Figure \ref{fig:thresholdregions}. Our program experienced significant slowdowns when either $a$ and $p$, or $b$ and $q$ were large, with the slowest occuring at the extreme point $a=b=p=q=0.975$. This suggests that a uniform bound when $\alpha\rightarrow1$ is impossible, and that Algorithm \ref{alg:optimalthreshold} will fail for sufficiently large values of $\alpha$. On the other hand, either $a=1$ or $b=1$ results in trivial limiting behaviour of the underlying chain, so this situation is unlikely to occur in any practical application.

Note that in Figure \ref{fig:thresholdregions}, the half of the main diagonal where $a+b>1$ is entirely blue. This can be proven analytically. We note that the condition $a+b>1$ corresponds to the underlying Markov chain having positive one-step autocorrelation, which is a reasonable assumption when the frequency of observation is greater than the frequency of change in the observed system, since in this case, one would not expect the system to oscillate with each observation.

\begin{proposition}
\label{prop:symmetric}
In the symmetric case $a=b$ and $p=q$, with $a+b>1$, the optimal general policy is the unique Region V policy $g_V$ given by $R_0\subset A_1$ and $R_1\subset A_0$.
\end{proposition}
\begin{proof}
Proposition \ref{prop:entropyformula} gives the limiting expected entropy as a convex combination of two quantities, hence
\begin{equation}
H_\infty\ge\min\Bigg\{\frac{\sum_kc_k^{(0)}H(z_k^{(0)})}{\sum_kc_k^{(0)}},\frac{\sum_kc_k^{(1)}H(z_k^{(1)})}{\sum_kc_k^{(1)}}\Bigg\}.
\label{eqn:symmetricH}
\end{equation}
Note that equality is realised when the two quantities above are equal, which occurs under $g_V$, since in that case $c_k^{(0)}=c_k^{(1)}$ and $z_k^{(0)}=1-z_k^{(1)}$.

Let $z_k=z_k^{(0)}$, $c_k=c_k^{(0)}/\sum c_k^{(0)}$ and $H_k=H(z_k^{(0)})$. In order to prove that $g_V$ is optimal, by symmetry, it suffices to prove that it minimises $\sum_kc_kH_k$.

First, we show that for each $k$ and any other policy $g$, $H_k(g_V)\le H_k(g)$, where the notation $H_k(g)$ makes the dependence on the policy explicit.

By Lemma \ref{lem:r0bigger}, $z_k(g)=r^k(0)\ge r_1^k(0)=z_k(g_V)$. By Lemma \ref{lem:monotonic} and our assumption that $a+b>1$, iterations of $r_0$ and $r_1$ approach their limits monotonically, hence $z_k(g)<\eta_0$ and $z_k(g_V)<\eta_1$. These inequalities are illustrated in Figure \ref{fig:symmetric}.

Combining these inequalities gives $z_k(g_V)\le z_k(g)<\eta_0$. Entropy is concave on the interval $[0,1]$, and therefore on the subinterval $[z_k(g_V),\eta_0]$, hence
\begin{equation}
\label{eqn:symmetricineq}
H(z_k(g))\ge\min\big\{H(z_k(g_V)),H(\eta_0)\big\}.
\end{equation}
Since $z_k(g_V)<\eta_1<\frac12$ and entropy is increasing on $[0,\frac12]$, $H(z_k(g_V))<H(\eta_1)$, which is equal to $H(\eta_0)$ by symmetry. Hence, the inequality (\ref{eqn:symmetricineq}) reduces to $H(z_k(g))\ge H(z_k(g_V))$, that is, $H_k(g)\ge H_k(g_V)$.

\vspace{3\parskip}
\begin{figure}[H]
  \centering
  \includegraphics[width=\textwidth]{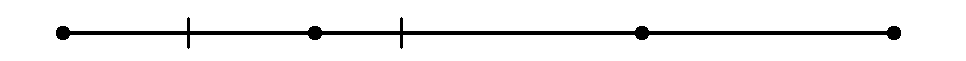}\\[-.5cm]
  0~~~~~~~~~~$z_k(g_V)$~~~~~~~$\eta_1$~~~~~$z_k(g)$ ~~~~~~~~~~~~~~~~~~~~$\eta_0$~~~~~~~~~~~~~~~~~~~~~~~~1~
  \large
  \caption[Proof of Proposition \ref{prop:symmetric}]{Diagram showing $z_k(g_V)$ and $z_k(g)$ in relation to 0, $\eta_1$, $\eta_0$ and 1. All positions are fixed except that $z_k(g)$ may be to the left of $\eta_1$. Since $z_k(g_V)$ lies to the left of $\eta_1$ and the diagram is symmetric, $z_k(g_V)$ has lower entropy than $\eta_0$. Since $z_k(g)$ lies between $z_k(g_V)$ and $\eta_0$ and entropy is concave, $z_k(g)$ has higher entropy than $z_k(g_V)$. Hence the $g_V$ minimises the entropy at $z_k$.}
  \label{fig:symmetric}
\end{figure}
\vspace\parskip

Next, we show that for each $k$ and any other policy $g$,
\begin{equation}
c_0(g_V)+c_1(g_V)+\cdots+c_k(g_V)\ge c_0(g)+c_1(g)+\cdots+c_k(g).
\label{eqn:coeffdominance}
\end{equation}
Let $a_k=\alpha(z_k^{(0)})$. Since $a+b>1$, by (\ref{eqn:specialcase}), $\alpha_1$ is decreasing, while $\alpha_0(z)=\alpha_1(1-z)$ is increasing. We have already established that $z_k(g_V)\le z_k(g)<\eta_0$ and $z_k(g_V)<\eta_1=1-\eta_0<1-z_k(g)$, which implies $\alpha_1(z_k(g_V))\ge\alpha_1(z_k(g))$ and $\alpha_1(z_k(g_V))\ge\alpha_0(z_k(g))$ respectively. This shows that $a_k(g_V)\ge a_k(g)$.

For $\ell<k$, write
\[c_0+\cdots+c_k=\frac{1+a_0+a_0a_1+\cdots+a_0a_1\cdots a_{k-1}}{1+a_0+a_0a_1+\cdots+a_0a_1\cdots a_{k-1}+\cdots}
\equiv\frac{X+\alpha_\ell Y}{X+\alpha_\ell Z}.\]
Since $X>0$ and $0<Y<Z$, decreasing $\alpha_\ell$ increases $c_0+\cdots+c_k$, and the same is true for $\ell\ge k$, since in that case we can write the expression in the same way with $Y=0$. This proves (\ref{eqn:coeffdominance}).

Using (\ref{eqn:coeffdominance}) and $H_0(g_V)<H_1(g_V)<H_2(g_V)<\cdots$, for any $\ell\in\N$,
\begin{align*}
S_\ell&\equiv\sum_{k\le\ell}\big(c_k(g_V)-c_k(g)\big)H_\ell(g_V)+\sum_{k>\ell}\big(c_k(g_V)-c_k(g)\big)H_k(g_V)\\
&\le\sum_{k\le\ell}\big(c_k(g_V)-c_k(g)\big)H_{\ell+1}(g_V)+\sum_{k>\ell}\big(c_k(g_V)-c_k(g)\big)H_k(g_V)=S_{\ell+1}.
\end{align*}
Since $\sum_kc_k=1$ identically, the second series vanishes as $\ell\rightarrow\infty$, while the first series is always non-negative by (\ref{eqn:coeffdominance}), hence $S_0\le0$. Thus,
\[\sum_kc_k(g_V)H_k(g_V)\le\sum_kc_k(g)H_k(g_V)\le\sum_kc_k(g)H_k(g).\]
This proves the required minimisation.
\end{proof}

Note that the proof above relies heavily on the fact that equality is attained in (\ref{eqn:symmetricH}). This occurs only in the symmetric case, and thus this approach does not generalise readily. The complexity of the proof in the symmetric case is indicative of the difficulty of the problem in general, and thus highlights the importance of the empirical description provided by Figure \ref{fig:thresholdregions}.

In the course of performing the computations to generate Figure \ref{fig:thresholdregions}, we noticed that entropy is unimodal with respect to threshold, with threshold considered as a circle as in Figure \ref{fig:thresholddomain}. While we cannot prove this analytically, it is true for each of the 152000 points in the parameter space considered.

This allows some simplification in finding the optimal threshold policy, since finding a local minimum is sufficient. Thus, we can alter Algorithm \ref{alg:optimalthreshold} to begin by testing only two policies, then testing policies in the direction of entropy decrease until a local minimum is found. However, the running time improvement is only a constant factor; if we model entropy as a standard sinusoid with respect to threshold, then the running time decreases by a factor of 3 on average.

\newpage
\subsection{General Policies}
The problem of determining the optimal general policy is much more difficult, due to the complexity of the space of general policies. Since a policy is uniquely determined by the value of the policy function at the orbit points, this space can be viewed as a hypercube of countably infinite dimension, which is much more difficult to study than the space of threshold policies, which is a circle.

One strategy is to truncate the orbit and consider a finite dimensional hypercube, justified by the fact that orbit points have masses which decay geometrically, and thus the tail contributes very little. However, a truncation at $N$ (that is, force the policy to be constant on $\big\{r^N(0),r^{N+1}(0),\ldots\big\}$, and similarly for the orbit of 1) gives $2^{2N+2}$ possible policies, which is still far too large to determine optimality by checking the entropy of each policy.

The next approximation is to only look for locally optimal policies, in the sense that changing the policy at each of the $2N+2$ truncated orbit points increases entropy, and hope that by finding enough such locally optimal policies, the globally optimal policy will be among them. Since a hypercube has very high connectivity, regions of attraction tend to be large, which heuristically suggests that this strategy will be effective.

\begin{algorithm}
\label{alg:localoptimum}
Finding a locally optimal truncated policy.
\vspace{-\parskip}
\begin{enumerate}
\item Pick $N$, and a starting policy, expressed as a pair of sequences of binary digits $g_k^{(i)}=g(z_k^{(i)})$, with $k=0,1,\ldots,N$;
\item Cycle through the digits $g_k^{(i)}$, flipping the digit if it gives a policy with lower entropy, otherwise leaving it unchanged;
\item If the previous step required any changes, repeat it, otherwise a locally optimal truncated policy has been found.
\end{enumerate}
\end{algorithm}

We picked $N=63$ since this allows a policy to be easily expressed as two unsigned 64-bit integers, and for each of the 152000 uniformly spaced parameters of Figure \ref{fig:thresholdregions}, we generated 10 policies uniformly on the hypercube and applied Algorithm \ref{alg:localoptimum}.

None of the locally optimal policies for any of the parameter values had lower entropy than the optimal threshold policy from Figure \ref{fig:thresholdregions}, and on average 98.3\% of them were equivalent to the optimal threshold policy, up to a precision of 0.1\%, indicating that the optimal threshold policy is locally optimal with a very large basin of attraction, which strongly suggests that it is also the globally optimal policy.

\begin{conjecture}
In the special case, the infimum of entropy attainable under threshold policies is the same as that under general policies.
\end{conjecture}

The fact that a large proportion of locally optimal policies have globally optimal entropy gives a new method for finding the optimal policy. By picking 10 random truncated policies and running Algorithm \ref{alg:localoptimum}, at least one of them will yield an optimal policy with very high probability. Empirical observations suggest that this method is slower than Algorithm \ref{alg:optimalthreshold} on average, but since the success rate remains high while Algorithm \ref{alg:optimalthreshold} becomes significantly slower as $\alpha$ approaches 1, this method is a better alternative for some parameter values.

\newpage
\vspace{3\parskip}
\begin{figure}[H]
  \centering
  \includegraphics[width=\textwidth]{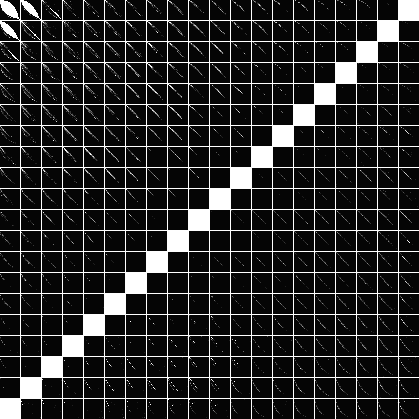}
  \large
  \caption[Locally optimal policies]{Locally optimal policies. Axes are as in Figure \ref{fig:thresholdregions}. Darkness increases with the proportion of simulated locally optimal policies which have the same entropy as the optimal threshold policy, up to a precision of 0.1\%. The average is 9.83 out of 10, but the distribution is far from uniform---local optima are exceedingly likely to be the same as the threshold optimum for some parameter values and exceedingly unlikely for others. The boundaries are approximately those of the threshold regions (see Figure \ref{fig:thresholdregions}), with some imprecision due to the non-deterministic nature of the simulation data.}
  \label{fig:localoptimum}
\end{figure}
\vspace\parskip
\newpage

One last policy of interest is the greedy policy. In the previous sections, we considered a long term optimality criterion in the minisation of limiting expected entropy, but in some cases, it may be more appropriate to set a short term goal. In particular, one may instead desire to minimise expected entropy at the next step, in an attempt to realise maximal immediate gain while ignoring future implications.

\begin{definition}
\label{def_greedy}
The \textbf{greedy} policy is the policy such that the expected entropy after one observation is minimised. Up to an exchage of strict and non-strict inequalities, this is given by
\begin{align*}
&A_0=\{z\in[0,1]:\alpha_0(z)H(r_0(z))<\alpha_1(z)H(r_1(z))\},\\
&A_1=\{z\in[0,1]:\alpha_0(z)H(r_0(z))\ge\alpha_1(z)H(r_1(z))\}.
\end{align*}
\end{definition}

The greedy policy has the benefit of being extremely easy to use, as it only requires a comparison of two functions at the current information state. Since these functions are smooth, efficient numerical methods such as Newton-Raphson can be used to determine the intersection points, thus allowing the policy to be described by a small number of thresholds.

In fact, only one threshold is required, as computational results show the greedy policy always a threshold policy. Using the 152000 uniformly distributed data points from before, in each case the two functions defining the greedy policy crossed over at most once.

\begin{conjecture}
The greedy policy is always a threshold policy.
\end{conjecture}
\begin{proof}[Idea of Proof]
Note that $q\alpha_0(z)r_0(z)=r_1(z)\alpha1(z)$. It may appear at first glance that the factor of $q$ violates symmetry, but recall that $z$ maps to $1-z$ under relabelling.

Using this identity, the intersection points that define the greedy policy satisfy $G(r_0(z))=qG(r_1(z))$, where $G(z)=H(z)/z$. It is easy to see that $G$ is monotonic decreasing on $[0,1]$ with range $[0,\infty]$, hence $F(z)=G^{-1}(qG(z))$ is a well-defined one-parameter family of functions mapping $[0,1]$ to itself with fixed points at 0 and 1.

On the other hand, $f(y)=r_0(r_1^{-1}(y))=y/\big((1-pq)y+pq\big)$ is also a one-parameter family of functions mapping $[0,1]$ to itself with fixed points at 0 and 1. Since the range of $r_0$ is contained in $(0,1)$, we can discount the endpoints $y=0$ and $y=1$, hence it suffices to show that the equation $F(y)=f(y)$ has at most one solution for $y\in(0,1)$.

Convexity methods may help in this last step but we have not been able to complete the proof.
\end{proof}

Even when the greedy policy is not optimal, it is very close to optimal. Of the 152000 uniformly distributed data points in Figure \ref{fig:greedyoptimal} below, the greedy policy is non-optimal at only 6698 points, or 4.41\%, up to an error tolerance of $10^{-12}$. On average, the greedy policy has entropy 0.0155\% higher than the optimal threshold policy, with a maximum error of 5.15\% occuring at the sample point $a=0.975$, $b=0.925$, $p=0.675$ and $q=0.375$. Thus the greedy polices provides an alternative suboptimal policy which is very easy to calculate and very close to optimal.

\newpage
\vspace{3\parskip}
\begin{figure}[H]
  \centering
  \includegraphics[width=\textwidth]{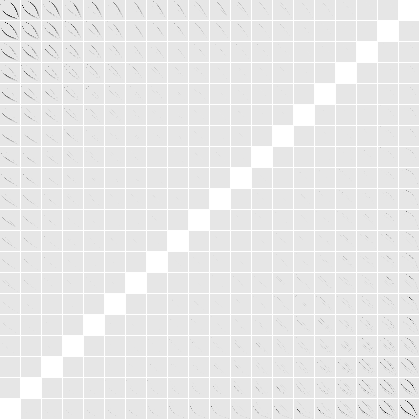}
  \large
  \caption[Optimality of the greedy policy]{Optimality of the greedy policy. Axes are as in Figure \ref{fig:thresholdregions}. Light grey indicates that the greedy policy is the optimal threshold policy; darker points indicate suboptimality with darkness proportional to error. Similarly to Figure \ref{fig:localoptimum}, the suboptimal points lie on the boundaries of the threshold regions.}
  \label{fig:greedyoptimal}
\end{figure}
\vspace\parskip

\newpage
We make a final remark that the likelihood of a locally optimal policy being globally optimal as shown in Figure \ref{fig:localoptimum}, and the closeness of the greedy policy to the optimal threshold policy as shown in Figure \ref{fig:greedyoptimal}, both exhibit a change in behaviour at the boundaries of the threshold regions as shown in Figure \ref{fig:thresholdregions}. This suggests that these regions are indeed qualitatively different, and are likely to be interesting objects for further study.

\newpage
\section{Conclusion and Future Work}
This thesis presents initial groundwork in the theory of hidden Markov models with multiple observation processes. We prove a condition under which the information state converges in distribution, and give algorithms for finding the optimal policy in a special case, which provides strong evidence that the optimal policy is a threshold policy.

Future work will aim to prove the conjectures that we were not able to prove:
\vspace{-2\parskip}
\begin{itemize}
\item The information state converges for an anchored hidden Markov model with only ergodicity rather than positivity;
\item The greedy policy is always a threshold policy;
\item Among threshold policies, the limiting expected entropy is unimodal with respect to threshold; and
\item The optimal threshold policy is also optimal among general policies.
\end{itemize}
Possible approaches to these problems are likely to be found in \cite{djonin}, \cite{legland} and \cite{macphee}. The author was not aware of these papers under after the date of submission, and thus was unable to incorporate their methods into this thesis.

Better algorithms and error bounds for finding the optimal policy are also a worthwhile goal. Although our algorithms are computationally feasible with reasonable prescribed errors, our focus was on finding workable rather than optimal algorithms, and thus there is plenty of room for improvement.

Another direction for future work would be to extend our results into more general cases.

\newpage
\addcontentsline{toc}{section}{Bibliography}
\bibliographystyle{plain}

\end{document}